\documentclass[11pt]{article}
\usepackage{amsmath}
\usepackage{amssymb}
\usepackage{enumitem}
\usepackage{theorem}
\usepackage[auth-sc-lg]{authblk}
\usepackage{graphicx} 
\input supp-pdf.tex 

\newtheorem{theorem}{Theorem}
\newtheorem{proposition}[theorem]{Proposition}
\newtheorem{corollary}[theorem]{Corollary}

\newcommand{\aut}{\mathrm{Aut}}
\newcommand{\aaut}{\mathrm{AAut}}
\newcommand{\caut}{\mathrm{CAut}}
\newcommand{\orb}{\mathrm{Orb}}
\newcommand{\rk}{\mathrm{rk}}
\newcommand{\AGL}{\mathrm{AGL}}

\title{Construction of infinite families of non-Schurian 
association schemes of order~$2p^2$, $p$ an odd prime, based on biaffine planes and
Heisenberg groups: research report and beyond}
\author[1]{\v Stefan Gy\"urki\footnote{e-mail addresses: {\tt gyurki@savbb.sk} (\v S. Gy\"urki), {\tt klin@cs.bgu.ac.il} (M. Klin)}}
\author[1,2]{Mikhail Klin}
\affil[1]{Matej Bel University} 
\affil[ ]{974 11 Bansk\'a Bystrica, Slovak Republic}
\affil[ ]{ }
\affil[2]{Department of Mathematics} 
\affil[ ]{Ben-Gurion University of the Negev}
\affil[ ]{84105 Beer Sheva, Israel}

\date{\today}

\begin{document}
\maketitle

\hrule

\begin{abstract}
Let $p$ be an odd prime. In this paper we provide a construction which gives
four non-Schurian association schemes for every $p\geq 5$ and two for $p=3$. This 
construction is explained using incidences between points and lines of a
biaffine plane and we also provide a pure algebraic model for it with the aid of 
finite Heisenberg groups.  
The obtained results are discussed in a more wide framework.
\medskip


\end{abstract}

\hrule

\section{Introduction}

This paper reports about research conducted on the edge between Algebraic
Graph Theory (briefly AGT) and Computer Algebra. 

Association schemes are one of the traditional areas of investigation in AGT. For a good decade, 
catalogues of small association schemes have been  available from
the web site \cite{hm}. It is known that all association schemes of order up
to 14 are Schurian, that is, they are coming from a sui\-table transitive 
permutation group in the standard manner. First examples of non-Schurian
association schemes exist on 15, 16 and 18 vertices. In particular, there are
just two classes of non-Schurian association schemes of order 18. 

D. Pasechnik explained in evident form in \cite{pa} how the non-Schurian rank 3
antisymmetric association scheme on 15 points appears,  and determined  its full
automorphism group.  

The famous Shrikhande graph generates non-Schurian rank 3 association schemes
on 16 points. In fact, there are many other non-Schurian schemes of order 16
(see \cite{hm}). Their clever computer free explanation might be of definite
interest. 

Surprisingly, such an explanation has not been given for the schemes on 18 points yet.
We are here filling this gap, providing an interpretation in terms of finite geometries.
Moreover, we are introducing a possible generalisation of these schemes which is leading
to four infinite families of non-Schurian association schemes on $2p^2$ points, where
$p>3$ is a prime. 
 
This paper originated from computer aided experiments, which were fulfilled by 
the author \v S.Gy.\ under the guidance of M.K. The manner in which the computer search 
helped to reach the presented results is, in our eyes, of independent interest.
The suggested structure of our paper was designed intentionally not only to reflect
the obtained scientific results, but also to pursue additional pedagogical, expository
and even philosophical goals.  

Preliminaries and the used computer tools are briefly discussed in 
Sections 2 and 3, respectively. Our starting geometrical object
is a biaffine coherent configuration $\mathcal M$, which appears via the 
intransitive permutation group $H\cong \mathbb Z_p^2\rtimes\mathbb Z_p$ 
of order $p^3$ and degree $2p^2$ having two orbits of length~$p^2$. Object $\mathcal M$ is
considered in Sections 4--8, along  with its four related color graphs~$\mathcal M_i$, 
$1\leq i\leq 4$, which turn out to be association schemes. In particular, 
an outline of a proof  that these association schemes are 
non-Schurian is given in Section 8. 

The presentation in Sections 4--8 is of a definite geometric nature, based on
a consideration of the classical biaffine plane of order $p$ with $p^2$ points
and $p^2$ lines.

In Section 9, all discussed structures are developed from
the scratch with the aid of an independent second model of $\mathcal M$, 
which is of a more algebraic nature. Here our starting structure is the finite
Heisenberg group of order $p^3$. 

Sections 10 and 11 are devoted to the use of algebraic groups of coherent 
configurations. Here we suggest some innovative methodological  elements, 
which might be of an independent interest for experts in AGT. 

Sections 12 and 13 are aimed to extend the scope of our  consideration, paying special 
attention to some extra features of our results, as well as to a number of well known
classes of graphs which appear to be relevant to new association schemes discovered 
by the authors. Some of the considered graphs play a significant role in Extremal Graph
Theory.

Section 14 serves as a research announcement of another portion of achieved results 
regarding new non-Schurian association schemes of constant rank 6 and 5, and some of their
properties. The level of rigour here is different from the main body of the text. In particular, 
no attempt is made to provide full formulations and outlines of proofs. 

Finally, the last Section 15 is written in a very specific style, quite 
typical of other publications co-authored by M.K. It is, in a sense, a mosaic of
diverse topics, implicitly or explicitly related to the content of the paper, however
not previously touched in evident form. Special attention is paid to the research paradigm 
``computer experiment, plausible reasonings, theoretical proof'', comparing views,
experience, and tastes of the authors with modern trends in science, which are correlated
with the exploitation of computer tools. 

Appendices 1--3 contain some routine data about the  investigated schemes. These data will play an essential role  in  
portions of the theoretical proofs.

\section{Preliminaries}

Below we provide a brief outline of the most significant concepts that will be 
used throughout the text. We refer to \cite{fk} and \cite{km} for a more
detailed background.

By a \emph{color graph} $\Gamma$ we will mean an ordered pair
$(V,\mathcal R)$, where $V$ is a set of vertices and~$\mathcal R$ 
a partition of $V\times V$ into binary relations. The elements of~$\mathcal R$ 
will be called \emph{colors}, and the number of
colors will be the \emph{rank} of $\Gamma$. In other words, a color graph
is an edge-colored complete directed graph with loops, whose arcs are colored
by the same color if and only if they belong to the same binary relation.  

A \emph{coherent configuration} is a color graph $\mathcal W=(\Omega,
\mathcal R)$, $\mathcal R=\{R_i \mid i\in I\}$, such that the following
axioms are satisfied:
\begin{itemize}
\item[(i)] The diagonal relation $\Delta_{\Omega}=\{(x,x)\mid x\in\Omega\}$
is a union of relations $\cup_{i\in I'}R_i$, for a suitable subset $I'
\subseteq I$. 
\item[(ii)] For each $i\in I$ there exists $i'\in I$ such that 
$R_i^T=R_{i'}$, where $R_i^T=\{(y,x)\mid (x,y)\in R_i\}$ is the relation 
transposed to $R_i$.
\item[(iii)] For any $i,j,k\in I$, the number $c_{i,j}^k$ of elements $z\in
\Omega$ such that $(x,z)\in R_i$ and $(z,y)\in R_j$ is a constant depending
only on $i,j,k$, and independent on the choice of $(x,y)\in R_k$. 
\end{itemize}
The numbers $c_{i,j}^k$ are called \emph{intersection numbers}, or 
sometimes \emph{structure constants} of $\mathcal W$. 

Assume that $|\Omega|=n$, and let us put $\Omega=\{1,2,\ldots,n\}$. To each \emph{basic graph}
$\Gamma_i=(\Omega, R_i)$ we associate its adjacency matrix $A_i=A(\Gamma_i)$. Then 
the set of \emph{basic matrices} $\{A_i\,|\, i\in I \}$ may be regarded as a basis of
a matrix algebra $\mathcal H$ which contains the identity matrix, the all-ones matrix~$J$, 
and is closed under transposition and Schur-Hadamard multiplication of matrices. 
Such an algebra is called a \emph{coherent algebra}, and we refer to the set $\{A_i\,|\, i\in I \}$
as its \emph{standard basis}.

The concepts of coherent configuration and coherent algebra were introduced 
by D. Higman (see e.g. \cite{hi}). Similar concepts were introduced independently by
B.Ju.\ Weisfeiler and A.A.\ Leman, see \cite{ws} and also \cite{rr} for an historical discussion.

A significant source of coherent configurations appears as follows. Assume that 
$(G,\Omega)$ is a permutation group acting on the set $\Omega$. For $(\alpha,\beta)\in
\Omega^2$ the set $\{(\alpha,\beta)^g\,|\, g\in G\}$, where $(\alpha,\beta)^g=
(\alpha^g,\beta^g)$, is called a \emph{2-orbit} of $G$, specifically the  \emph{2-orbit} of $G$ 
\emph{corresponding to} $(\alpha,\beta)$. (Note that when  $(G,\Omega)$ is a transitive permutation
group, many authors prefer the term \emph{orbital} for this set.)   

Denoting by 2-$\orb(G,\Omega)$ the set of 2-orbits of a permutation group $(G,
\Omega)$, it is easy to check that $(\Omega,\;$2-$\orb(\Omega))$ is a
coherent configuration. Coherent configurations that arise 
in this manner are called \emph{Schurian}, otherwise we call them \emph{non-Schurian}. 

An \emph{association scheme} $\mathcal W=(\Omega,\mathcal R)$ (also called a  
\emph{homogeneous coherent configuration}) is a coherent configuration in which the diagonal relation
$\Delta_{\Omega}$  belongs to $\mathcal R$. Thus, Schurian association
schemes are coming from transitive permutation groups.

A coherent configuration $\mathcal W$ is called \emph{commutative}  
if for all $i,j,k\in I$ we have $c_{ij}^k=c_{ji}^k$.  We call $\mathcal W$  
 \emph{symmetric} if $R_i=R_i^T$ for all
$i\in I$. It is a well known fact that a symmetric coherent configuration
is also commutative, but the converse is not true in general. 

To each coherent configuration $\mathcal W$ we may assign three groups:
$\aut(\mathcal W), \caut(\mathcal W)$ and $\aaut(\mathcal W)$. 
The (combinatorial) \emph{group of automorphisms} $\aut(\mathcal W)$
consists of the permutations $\phi:\Omega\to\Omega$ which preserve
the relations, i.e. $R_i^{\phi}=R_i$ for all $R_i\in\mathcal R$.
The \emph{color automorphisms} are permitted to permute the relations from $\mathcal R$,  i.e. 
for $\phi:\Omega\to\Omega$ we have $\phi\in\caut(\mathcal W)$ if and only if
for all $i\in I$ there exists $j\in I$ such that $R_i^{\phi}=R_j$. 
An \emph{algebraic automorphism} is a bijection $\psi:\mathcal R\to
\mathcal R$ that satisfies $c_{ij}^k=c_{i^{\psi}j^{\psi}}^{k^{\psi}}$.
It is easy to verify that $\aut(\mathcal W)$ is a normal subgroup of $\caut(\mathcal W)$, 
and that the quotient  group $\caut(\mathcal W)/\aut(\mathcal W)$ embeds 
naturally in $\aaut(\mathcal W)$. 

For each group $K$ of algebraic automorphisms of $\mathcal W=(\Omega,
\mathcal R)$ one can define an \emph{algebraic merging} of 
$\mathcal R$ in the following way. 
Let $\mathcal R/K$ denote the set of orbits of $K$ on
$\mathcal R$. For each $O\in \mathcal R/K$ define $O^+$ to be the union
of all relations from $O$. Then the set of relations $\{O^+\,|\, O\in
\mathcal R/K\}$ forms a coherent configuration on $\Omega$. We will call it an 
\emph{algebraic merging} of $\mathcal R$ with respect to $K$. Note that
if $K\leq \caut(\mathcal W)/\aut(\mathcal W)$ and  $\mathcal W$ is
Schurian, then the resulting merging is Schurian as well. In contrast,
a merging with respect to a subgroup of $K$ not contained in 
$\caut(\mathcal W)/\aut(\mathcal W)$ may lead to a non-Schurian
coherent configuration. 

It is clear from the definitions that an association scheme $\mathcal W$ is
Schurian if and only if its rank coincides with the rank of its group of
automorphisms $\aut(\mathcal W)$.

\section{Computer tools}

This project heavily depends on the use of computer tools. Moreover, in a
sense it can serve as a pattern for the use of computers in AGT. 
Indeed, we made a lot of fast computations; investigated and organized the obtained results;
analyzed all relevant data; made conjectures based on these data; 
checked these conjectures for higher values
of parameters of the considered series; transformed these  conjectures to a suitable analytical form;
transformed the obtained results from numerical to symbolic mode; 
created pictures and tables; and so on.
 
Here we are mainly working with coherent configurations and association schemes, as well as 
with the permutation groups related to them. For this purpose, in 1990--92 a computer package was
created in Moscow as a result of the activities of I.A.\ Farad\v zev and the author M.K.  
This package goes by the name {\sf COCO}, and was introduced in \cite{co}; 
see also \cite{fk} for deeper consideration of the used methodology and algorithms. 
{\sf COCO} is still very helpful for performing initial computational experiments.

Nevertheless, nowadays the mainstream of our computer aided activities is based on the
use of the free software {\sf GAP} \cite{ga} (Groups, Algorithms and Programming),   
in particular its share package {\sf GRAPE} \cite{gr} which   
works in conjunction with {\sf nauty} \cite{na}. 

In addition, we strongly benefitted from the kind permission  granted by Sven Reichard 
to use his unpublished package  {\sf COCO~IIR} \cite{sv}  (still under development) which operates under the {\sf GAP} platform. The foremost goal of {\sf COCO~IIR} is to extend the scope of 
algorithms from {\sf COCO}, relying on many new developments
and fresh ideas adopted from modern computer algebra. 

Last but not least, we are pleased to acknowledge a package of programs for
computing with association schemes, written by Hanaki and Miyamoto. It
is open-source software that may be freely downloaded  from the homepage of the authors
\cite{hm}. This package also works under {\sf GAP}, and contains some tools that have
not yet been implemented in {\sf COCO~IIR}.

\section{A biaffine coherent configuration from the biaffine plane}

Let $p$ be an odd prime, and let $\mathbb Z_p$ be the cyclic group of order $p$. 
Throughout this text, the set of nonzero elements in $\mathbb Z_p$ 
will be denoted by $\mathbb Z_p^*$. 
Take two copies $\mathcal P$ and $\mathcal L$ of $\mathbb Z_p\times \mathbb Z_p$. 
The first copy $\mathcal P$ is none other than  the point set of the classical 
(Desarguesian) affine plane of order $p$. 
Each element $P\in\mathcal P$ may be identified uniquely 
with a pair of coordinates of the form $P=[x,y]$.
 Thus, we refer to the elements of $\mathcal P$ as  \emph{points}. 
Let $\mathcal L$ be the set of \emph{"non-vertical" lines} in the affine plane, 
i.e. $\ell\in\mathcal L$ if and only if the equation for $\ell$ may be expressed as 
$y=k\cdot x + q$ for some $k,q\in \mathbb Z_p$. 
Each line $\ell$ is determined uniquely by a pair $\ell=(k,q)$.  
In order to distinguish points and lines we will use square brackets for
points and parentheses for lines. 

\begin{figure}[htb]
\begin{center}
\includegraphics[width=13cm]{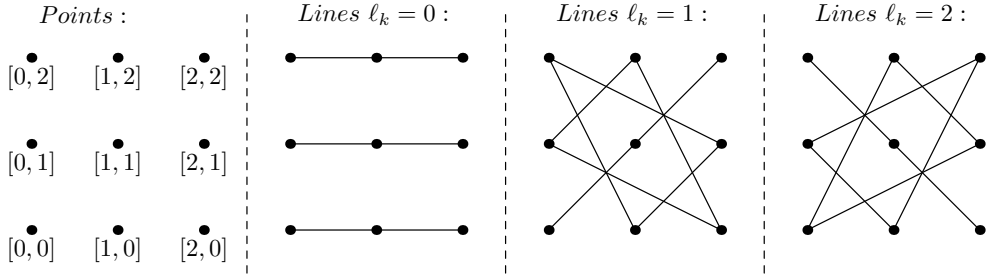}
\end{center}
\caption{The objects of the biaffine plane $\mathcal B_3$.}
\label{fig.1}
\end{figure}

This geometry, in which one parallel class of the affine plane has been removed 
(in our  case, the class of ``vertical'' lines) goes by the  
well established name \emph{biaffine plane}.
We denote it by~$\mathcal B_p$ (see Figure~\ref{fig.1} depicting~$\mathcal B_3$). 
For more details about biaffine planes and their relation to
other combinatorial structures, we recommend \cite{wi};  see also additional comments
at the end of this paper. 

Given a point $P=[x,y]$ and a line $\ell=(k,q)$ of the biaffine plane, 
we define a \emph{quasidistance} 
$d:(\mathcal P\times\mathcal L)\cup(\mathcal L\times\mathcal P)\to \mathbb Z_p$ 
by the formulas: $d(P,\ell)=k\cdot x+q-y$ and $d(\ell,P)=y-q-k\cdot x$. 
Note that $d$  does not define a metric (indeed, it is not symmetric 
and non-negativity does not make sense in $\mathbb Z_p$).  
Thus, it is just a vague analogue. 

Recall that a natural source of coherent configurations is coming from permutation
groups. If we take any permutation group $(G,\Omega)$, then $\mathcal W=(\Omega,
\;$ 2-$\orb(G))$ is its associated (Schurian) coherent configuration. 
In particular, if $(G,\Omega)$ is transitive, then $\mathcal W$
is an association scheme. 

Let us now consider an action of  the group $H=(\mathbb Z_p)^2\rtimes\mathbb Z_p$ 
on the set $\Omega=\mathcal P\cup\mathcal L$, most conveniently described in 
terms of generators. To each pair $(a,b)\in\mathbb Z_p^2$ we associate a 
\emph{translation}~$t_{ab}$ acting naturally on $\mathcal P$ as
$[x,y]\mapsto [x+a,y+b]$, while the induced action on $\mathcal L$ is
$(k,q) \mapsto (k, b+q-ak)$. Of course, the set of all translations forms an Abelian
group of order $p^2$ under composition of permutations (denoted $\circ$)  
and is isomorphic to $\mathbb Z_p^2$.
Further, let $\varphi: \mathcal P\to\mathcal P$ be defined by  
$\varphi: [x,y]\mapsto [x,y-x]$. The corresponding permutation on lines is
$(k,q)\mapsto (k-1,q)$. Clearly, $\varphi$ is a permutation of order $p$ and
it is immediate that $\varphi\circ t_{a,b-a}=t_{a,b}\circ \varphi$. 
Our group $H$ above is generated by all translations together with the permutation $\varphi$.
Note that all elements of $H$ may be expressed in the form $\varphi^u\circ t_{a,b}$,
where $a,b,u$ are suitable elements of~$\mathbb Z_p$. Moreover, to distinct 
triples of $a,b,u$ correspond distinct elements of $H$. 
In other words, $H=\langle t_{a,b}, \varphi\rangle$ and $|H|=p^3$.
Now set $h_{a,b,u}:=\varphi^u\circ t_{a,b}$. Then  multiplication in $H$
is given by $h_{a,b,u}\circ h_{c,d,v}=h_{a+c, b+d-av, u+v}$. 

Observe that the action of $H$ is intransitive on $\Omega$ with two orbits 
$\mathcal P$ and $\mathcal L$. 

In what follows, the introduced permutation group $(H,\Omega)$ will be called the 
{\it Heisenberg group modulo $p$}; see Section 9 for more details. 

\begin{proposition}\label{prop1}
Let $p$ be an odd prime. Then 
\begin{itemize}
\item[(1)] the rank of $(H,\Omega)$ is $6p-2$,
\item[(2)] the $2$-orbits of $(H,\Omega)$ may be divided into six different 
types of classes $A_i, B_i, C_i, D_i, E_i$ and $F_i$, which are characterized by
suitable relations between coordinates of objects in $\Omega$. 
\end{itemize}
\end{proposition}

\noindent{\bf Proof.} 
Let $P_1=[x_1,y_1]$, $P_2=[x_2,y_2]\in\mathcal P$ and $\ell_1=(k_1,q_1)$, 
$\ell_2=(k_2,q_2)\in\mathcal L$. Then the types of classes are the following:
\begin{itemize}
\item $(P_1,P_2)\in A_i \iff x_1=x_2$ and $y_2-y_1=i$, where $i\in\mathbb Z_p$
(i.e., those pairs of points whose first coordinates are equal and second coordinates differ
by $i$ in the given order; note that $A_0$ is the diagonal relation on $\mathcal P$),
\item $(P_1,P_2)\in B_i \iff x_2-x_1=i$, where $i\in\mathbb Z_p^*$
(i.e., those pairs of points whose first coordinates differ by $i$ in the given order),
\item $(\ell_1,\ell_2)\in C_i \iff k_1=k_2$ and $q_2-q_1=i$, where $i\in\mathbb Z_p$ 
(i.e., those pairs of lines whose first coordinates are equal and second coordinates differ
by $i$ in the given order; note that $C_0$ is the diagonal relation on $\mathcal L$),
\item $(\ell_1,\ell_2)\in D_i \iff k_2-k_1=i$, where $i\in\mathbb Z_p^*$
(i.e., those pairs of lines whose first coordinates differ by $i$ in the given order),
\item $(P_1,\ell_1)\in E_i \iff k_1\cdot x_1+q_1-y_1=i$, where $i\in\mathbb Z_p$
(i.e., those point-line pairs whose quasidistance $d(P_1,\ell_1)$ is $i$),
\item $(\ell_1,P_1)\in F_i \iff y_1-k_1x_1-q_1=i$, where $i\in\mathbb Z_p$
(i.e., those point-line pairs whose quasidistance $d(\ell_1,P_1)$ is $i$).
\end{itemize}

To complete the proof, one must verify two things: that the $6p-2$ relations 
introduced above  indeed form a partition of the set $\Omega^2$, and that each 
such relation is in fact a 2-orbit of $(H,\Omega)$. 

\rightline{$\square$}

\noindent{\bf Remark 1.} One can easily check from its definition that the permutation 
$\varphi$ has exactly $p$ fixed points in its action on $\mathcal P$, as well as $p$
fixed points in its action on $\mathcal L$. Relying on the bijection between 2-orbits 
of a transitive permutation group and orbits (1-orbits) of the stabilizer of an 
arbitrary point (see e.g. \cite{fk}), the reader can easily deduce that there exists 
exactly $p+(p-1)$ 2-orbits of the transitive action $(H,\mathcal P)$, and similarly
for the action $(H,\mathcal L)$. Thus, we obtain in this manner $4p-2$ 2-orbits of 
$(H,\Omega)$. Observing that there are $p$ 2-orbits of type $E_i$ and $p$ of type $F_i$, 
we arrive at the desired amount of $6p-2$. 

\medskip

\noindent{\bf Remark 2.} Reflexive 2-orbits $A_0$ and $C_0$ are obviously symmetric, 
however all remaining 2-orbits are antisymmetric. Namely, we obtain that 
$A_i^T=A_{p-i}$, $B_i^T=B_{p-i}$, $C_i^T=C_{p-i}$, $D_i^T=D_{p-i}$, $E_i^T=F_{p-i}$ 
and $F_i^T=E_{p-i}$. Here and below, operations on subscripts are performed as in 
$\mathbb Z_p$. 

This description provides a quite nice geometrical interpretation of the relations. 
In the sections to follow, we will define certain families of color graphs as 
suitable mergings of these relations. We will prove that these color graphs are 
in fact association schemes, and we will investigate them further.
  
Note that at the end of the paper we will once again justify a portion of our 
obtained results, this time relying  on a more  purely algebraic approach.

Denote by $\mathcal M$ the coherent configuration corresponding to the group $H$ 
in our construction above. For obvious reasons, we shall refer to $\mathcal M$ as 
a \emph{biaffine coherent configuration}.
  
\section{Intersection numbers of the biaffine coherent configuration $\mathcal M$}

The biaffine coherent configuration $\mathcal M$ was constructed with the
aid of a group, that is, it is Schurian.  Presently, we are interested in 
the intersection numbers of $\mathcal M$. 
In this section we display these numbers with the aid of  tables.
In each individual table the superscript is fixed, the symbol in the row indicates
the first subscript, and the symbol in the column indicates the second subscript. 


We are using the \emph{Kronecker's symbol} $\delta_{i,j}$
in order to shorten computations and formulas. 

\begin{proposition}
The tensor of structure constants of the biaffine coherent configuration
$\mathcal M$ is given as follows:
\vspace*{4pt}

\begin{minipage}{7cm}
\begin{tabular}{|c|ccc|}
\hline
$c_{r_i,c_j}^{A_k}$ & $A_j$ & $B_j$ & $F_j$ \\
\hline
$A_i$ & $\delta_{i+j,k}$ & $0$ & $0$ \\
$B_i$ & $0$ & $p\cdot\delta_{i+j,0}$ & $0$ \\
$E_i$ & $0$ & $0$ & $p\cdot\delta_{i+j,k}$ \\
\hline
\end{tabular}
\end{minipage}
\begin{minipage}{7cm}
\begin{tabular}{|c|ccc|}
\hline
$c_{r_i,c_j}^{C_k}$ & $C_j$ & $D_j$ & $E_j$ \\
\hline
$C_i$ & $\delta_{i+j,k}$ & $0$ & $0$ \\
$D_i$ & $0$ & $p\cdot\delta_{i+j,0}$ & $0$ \\
$F_i$ & $0$ & $0$ & $p\cdot\delta_{i+j,k}$ \\
\hline
\end{tabular}
\end{minipage}

\medskip

\begin{minipage}{7cm}
\begin{tabular}{|c|ccc|}
\hline
$c_{r_i,c_j}^{B_k}$ & $A_j$ & $B_j$ & $F_j$ \\
\hline
$A_i$ & $0$ & $\delta_{j,k}$ & $0$ \\
$B_i$ & $\delta_{i,k}$ & $p\cdot\delta_{i+j,k}$ & $0$ \\
$E_i$ & $0$ & $0$ & $1$ \\
\hline
\end{tabular}
\end{minipage}
\begin{minipage}{7cm}
\begin{tabular}{|c|ccc|}
\hline
$c_{r_i,c_j}^{D_k}$ & $C_j$ & $D_j$ & $E_j$ \\
\hline
$C_i$ & $0$ & $\delta_{j,k}$ & $0$ \\
$D_i$ & $\delta_{i,k}$ & $p\cdot\delta_{i+j,k}$ & $0$ \\
$F_i$ & $0$ & $0$ & $1$ \\
\hline
\end{tabular}
\end{minipage}

\medskip

\begin{minipage}{7cm}
\begin{tabular}{|c|ccc|}
\hline
$c_{r_i,c_j}^{E_k}$ & $C_j$ & $D_j$ & $E_j$ \\
\hline
$A_i$ & $0$ & $0$ & $\delta_{i+j,k}$ \\
$B_i$ & $0$ & $0$ & $1$ \\
$E_i$ & $p\cdot\delta_{i+j,k}$ & $1$ & $0$ \\
\hline
\end{tabular}
\end{minipage}
\begin{minipage}{7cm}
\begin{tabular}{|c|ccc|}
\hline
$c_{r_i,c_j}^{F_k}$ & $A_j$ & $B_j$ & $F_j$ \\
\hline
$C_i$ & $0$ & $0$ & $\delta_{i+j,k}$ \\
$D_i$ & $0$ & $0$ & $1$ \\
$F_i$ & $p\cdot\delta_{i+j,k}$ & $1$ & $0$ \\
\hline
\end{tabular}
\end{minipage}

\medskip

\noindent where the indices $i,j,k$ go through all feasible values. 
All structure constants not displayed here  are zero.

\end{proposition}

\noindent {\bf Example 1.} 
In the biaffine coherent configuration of order 50 (i.e. $p=5$)
we have 
\begin{align*}
c_{A_2,A_4}^{A_1}&=\delta_{2+4,1}=\delta_{1,1}=1, \\
c_{B_2,B_3}^{A_4}&=5\cdot\delta_{2+3,0}=5\cdot\delta_{0,0}=5, \\
c_{C_2,A_4}^{A_3}&=0, \\
c_{D_2,F_0}^{F_1}&=1. 
\end{align*}

\medskip

\noindent{\bf Proof (outline).}

First observe that for all $i\in \mathbb Z_p$ and $j\in\mathbb Z_p^{\ast}$,  
we have $A_i,B_j\subseteq \mathcal P\times\mathcal P$,
$C_i,D_j\subseteq \mathcal L\times\mathcal L$,
$E_i\subseteq \mathcal P\times\mathcal L$ and
$F_i\subseteq \mathcal L\times\mathcal P$. 
 Thus  all structure constants of the form $c_{Y_i,Z_j}^{X_k}$ are zero provided $X,Y,Z$  satisfy any of the following:
\begin{itemize}
\item $X\in\{A,B,E\}$ and $Y\in\{C,D,F\}$, or  $Y\in\{A,B,E\}$ and $X\in\{C,D,F\}$,
\item $X\in\{A,B,F\}$ and $Z\in\{C,D,E\}$, or  $Z\in\{A,B,F\}$ and $X\in\{C,D,E\}$,
\item $Y\in\{A,B,F\}$ and $Z\in\{C,D,F\}$, or  $Y\in\{C,D,E\}$ and $Z\in\{A,B,E\}$.
\end{itemize}
The reason is simply that in these cases the composition of relations $Y_i$ and $Z_j$ 
is either impossible, or giving a relation disjoint to $X_k$. This observation is crucial
in order to understand that those structure constants  which in principle may be nonzero 
are covered just by the six kinds of tables presented above.

Simple algebraic manipulations of coordinates lead us to the following:
$$c_{A_i,A_j}^{A_k}=c_{C_i,C_j}^{C_k}=\delta_{i+j,k}, \quad c_{B_i,B_j}^{A_k}=c_{D_i,D_j}^{C_k}=
p\cdot\delta_{i+j,0}, \quad c_{A_i,B_j}^{B_k}=c_{C_i,D_j}^{D_k}=\delta_{j,k},$$ 
$$c_{B_i,A_j}^{B_k}=c_{D_i,C_j}^{D_k}=\delta_{i,k}, \textrm{ and } 
c_{B_i,B_j}^{B_k}=c_{D_i,D_j}^{D_k}=p\cdot\delta_{i+j,k}.$$

Computation of the remaining structure constants is not as  straightforward and 
requires a bit more sophistication. Yet, these may be determined by counting 
with coordinates. For example, to show that $c_{E_i,F_j}^{A_k}=p\cdot\delta_{i+j,k}$
we first consider a pair $(P_1,P_2)\in A_k$. Putting $P_1=[x_1,y_1]$, this means 
that $P_2$ is uniquely determined by the coordinates $P_2=[x_1,y_1+k]$. We now seek 
the number of lines $\ell=(m,q)$ for which $(P_1,\ell)\in E_i$ and $(\ell,P_2)\in F_j$. 
This yields two equations that must be satisfied by the coordinates of $\ell$:
\begin{align*}
m\cdot x_1 + q - y_1 & = i \\
(y_1+k) - m\cdot x_1 - q &= j.
\end{align*} 

\noindent The sum of these two equations tells us that there are no solutions
when $i+j\neq k$ in $\mathbb Z_p$. To the contrary, if $i+j=k$ in $\mathbb Z_p$ 
then all lines that satisfy the first equation are  automatically 
solutions to the entire system. In this manner we obtain precisely $p$ solutions, one  
for each fixed choice of $m\in\mathbb Z_p$. 

In a similar fashion one can derive all remaining structure constants of $\mathcal M$. 

\rightline{$\square$}

\section{Construction of four families of color graphs}

Now we are ready to describe our four families of color graphs. 
To emphasize the connection between these color graphs and the initial 
biaffine coherent configuration $\mathcal M$, we define basic relations 
in terms of  $A_i,B_i,C_i,D_i,E_i$ and $F_i$. Later we will show that these 
color graphs are corresponding to association schemes. 

Let us consider the following subsets of $\Omega\times\Omega$:

\begin{itemize}
\item $R_0=A_0\cup C_0$,
\item $S_i=A_i\cup C_i$, where $i\in\mathbb Z_p^*$, 
\item $T_i=B_i\cup D_i$, where $i\in\mathbb Z_p^*$, 
\item $U_i=E_i\cup F_i$, where $i\in\mathbb Z_p$.
Note that the relation $U_0$ coincides with the set of  flags in the biaffine plane
$\mathcal B_p$. 
\end{itemize}

Further, let $S_i^*=S_i\cup S_{p-i}$, $T_i^*=T_i\cup T_{p-i}$ and $U_i^*=U_i
\cup U_{p-i}$ be the respective symmetrizations of the relations $S_i$, $T_i$ and $U_i$, 
canonically denoted for each $i\in\{1,2,\ldots,(p-1)/2\}$. Finally, 
let $S=S_1\cup S_2\cup \ldots\cup S_{p-1}$ and  $U=U_1\cup U_2\cup\ldots\cup 
U_{p-1}$.  Observe that $U$ is the set of antiflags in $\mathcal B_p$. 

It is straightforward to check that $\{R_0,S_1,\ldots,S_{p-1},T_1,\ldots,T_{p-1},U_0,
U_1,\ldots,U_{p-1}\}$ forms a partition of $\Omega\times\Omega$.

It remains to define the requested color graphs $\mathcal M_i$, $1\leq i\leq 4$,  
each with vertex set $\Omega$.

\noindent{\bf Color graph 1.} Denote by $\mathcal M_1$  the color graph with 
colors given by the sets $R_0$, $S_1,\ldots,S_{p-1}$, $T_1,\ldots,T_{p-1}$,
$U_0$, $U_1,\ldots,U_{p-1}$. 
 
\noindent{\bf Color graph 2.} Denote by $\mathcal M_2$ the color graph with  
colors given by $R_0$, $S_1^*$, $S_2^*, \ldots, S_{(p-1)/2}^*$, $T_1$, $T_2, \ldots, 
T_{p-1}$, $U_0$, $U_1^*$, $U_2^*, \ldots, U_{(p-1)/2}^*$. 

\noindent{\bf Color graph 3.} Denote by $\mathcal M_3$ the color graph with  
colors given by  $R_0$, $S$, $T_1$, $T_2,\ldots,T_{p-1}$, $U_0$, $U$. 

\noindent{\bf Color graph 4.} Finally, denote by $\mathcal M_4$ the color 
graph with colors given by $R_0$, $S$, $T_1^*$, $T_2^*, \ldots, T_{(p-1)/2}^*$, 
$U_0$, $U$. 

\smallskip

Note that for $p=3$ the color graphs $\mathcal M_2$ and $\mathcal M_3$
coincide. 

These four color graphs will play an important role in this paper. 

\section{Intersection numbers}

We wish now to show that our color graphs $\mathcal M_1$--$\mathcal M_4$ 
defined in the previous section are association schemes. 
Observe that it is sufficient to show the existence of intersection numbers 
(structure constants) because all other axioms of an  association scheme are
trivially satisfied. 

For convenience, we shall invoke the following notational simplification for 
the indices of intersection numbers. We shall write  $si, ti, ui$ in place of 
$S_i$, $T_i$, $U_i$, respectively. For example, $c_{si,tj}^{uk}$ indicates the 
number of elements $z\in\Omega$ such that $(x,z)\in S_i$ and $(z,y)\in T_j$ 
for any $(x,y)\in U_k$. A subscripted or superscripted zero shall always 
indicate  the relation $R_0$, while any index $i$ not accompanied by a 
specified symbol will indicate any feasible relation. For the sake of brevity, 
we shall only indicate those intersection numbers in which starred relations 
(such as $U_k^*$) occur  in the superscript. In each case we are providing 
only one argumentation (usually for points), because the dual consideration 
(for lines) is similar. We make frequent use of \emph{Kronecker's symbol} 
$\delta_{i,j}$ in order to shorten computations and formulas. 

To make enumeration easier the following observations are helpful:

\medskip

\noindent{\bf Observation 1.} For all $1\leq i\leq p-1$ and $1\leq j\leq 
(p-1)/2$ we have $R_0,S,S_i,T_i,S_j^*,T_j^*\subseteq (\mathcal P\times
\mathcal P)\cup(\mathcal L\times\mathcal L)$, and $U,U_0,U_i,U_j^*\subseteq
(\mathcal P\times\mathcal L)\cup(\mathcal L\times\mathcal P)$. 
Thus, the intersection numbers of type 
$c_{si,sj}^{uk},c_{si,tj}^{uk},c_{ti,sj}^{uk},c_{ti,tj}^{uk},c_{ui,uj}^{uk},
c_{si,uj}^{sk},c_{si,uj}^{tk},c_{ti,uj}^{sk},c_{ti,uj}^{tk},
c_{ui,sj}^{sk},c_{ui,sj}^{tk},c_{ui,tj}^{sk},c_{ui,tj}^{tk}$ are zero 
for all choices of  $i,j,k$. 

\medskip

\noindent{\bf Remark 3.} Let us choose the symbol $\ast$ to indicate \emph{composition}
of relations. Recall that in a coherent algebra this operation corresponds to a
product of corresponding adjacency matrices. In a coherent configuration the result
of composition is usually a \emph{multirelation}, that is, a set of relations together
with their (non-negative) integer multiplicities. Thus we have adopted $\ast$ in 
order to avoid misunderstandings, since the binary operator $\circ$ is reserved 
for the Schur-Hadamard product in the theory of coherent configurations.

\medskip

\noindent{\bf Observation 2.} For the compositions of relations 
$S_i,S_j,T_i,T_j$ we have $S_i\ast T_j=T_j$, $T_i\ast S_j=T_i$ and if 
$i+j\neq 0$, then $S_i\ast S_j=S_{i+j}$ and $T_i\ast T_j=T_{i+j}$. 
As a consequence we obtain the following: $c_{si,tj}^{sk}=c_{tj,si}^{sk}=
c_{si,sj}^{tk}=0$, and for $i+j\neq 0$: $c_{ti,tj}^{sk}=0$, 
$c_{si,sj}^{sk}=\delta_{i+j,k}$.

\medskip

\noindent{\bf Observation 3.} 
For each color $X$ we have $R_0\ast X=X\ast R_0=X$, and for any 
$Y\neq X^T$ we obtain $X\ast Y\neq R_0\neq Y\ast X$. Thus for $i\neq j$: 
$c_{0,i}^j=c_{i,0}^j=0$, and for $i\neq j'$: $c^0_{i,j}=0$. 

\medskip

\noindent{\bf Observation 4.} Let $P_1,P_2\in\mathcal P$ be two collinear 
points in $\mathcal B_p$, and let $L_1=\{\ell\in\mathcal L\,|, P_1\in\ell\}$.
Then for all $\ell_i,\ell_j\in L_1$ we have $d(P_2,\ell_i)=d(P_2,\ell_j)$ if
and only if $i=j$, where $d(P,l)$ is the previously defined quasidistance.

\medskip

All intersection numbers are displayed in Appendix 1. These values were 
derived by geometrical arguments, usually by considering points and lines at 
a given quasidistance from two objects. For example, consider the color graph
$\mathcal M_1$, and let $P=[x,y]$ and $i\in \mathbb Z_p$ be fixed. Then there 
is a unique point with coordinates $[x,y+i]$, precisely $p$ points with first 
coordinate $x+i$, and precisely $p$ lines at quasidistance $i$ from $P$.  
These observations lead directly to the intersection numbers 
$c_{si,sj}^0=\delta_{i,j'}$,  $c_{ti,tj}^0=p\cdot\delta_{i,j'}$ and   
$c_{ui,uj}^0=p\cdot\delta_{i,j'}$, respectively.

\medskip

\begin{theorem}\label{thm3}
The following holds:
\begin{itemize}
\item[(a)]  $\mathcal M_1,\mathcal M_2,\mathcal M_3,\mathcal M_4$ 
are association schemes.
\item[(b)] The combinatorial groups of automorphisms of 
$\mathcal M_1,\mathcal M_2,\mathcal M_3,\mathcal M_4$ contain 
a subgroup isomorphic to $H$.
\item[(c)] $\mathcal M_2$ is a merging of $\mathcal M_1$,
$\mathcal M_3$ is a merging of $\mathcal M_2$, and 
$\mathcal M_4$ is a merging of $\mathcal M_3$. 
\item[(d)] 
\[
\begin{array}{rcl}
\rk(\mathcal M_1)&=&3p-1, \\ 
\rk(\mathcal M_2)&=&2p, \\ 
\rk(\mathcal M_3)&=&p+3, \\
\rk(\mathcal M_4)&=&(p+7)/2.\\
\end{array}
\]
\end{itemize}
\end{theorem}

\noindent{\bf Proof.} 
Parts (a) and (b) have already been proven. Proofs of (c) and (d) follow easily 
from the definition of $\mathcal M_i$, $1\le i\le 4$.
  
\rightline{$\square$}

\section{Automorphism groups of the association schemes}

Recall that to each association scheme $\mathcal M$ we may assign three 
groups: $\aut(\mathcal M), \caut(\mathcal M)$ and $\aaut(\mathcal M)$. 
In this section we will focus on the combinatorial group of automorphisms~$\aut(\mathcal M)$. 
Recall that this group consists of all permutations $\phi:\Omega\to\Omega$ that preserve
relations, i.e. $R_i^{\phi}=R_i$ for all $R_i\in\mathcal R$.

\begin{theorem} \label{thm4}
Let $\aut(\mathcal M_1)$, $\aut(\mathcal M_2)$, $\aut(\mathcal M_3)$ and 
$\aut(\mathcal M_4)$ be the combinatorial groups of automorphisms of $\mathcal M_1$, 
$\mathcal M_2$, $\mathcal M_3$ and $\mathcal M_4$, respectively. Then the following hold:
\begin{itemize}
\item[(a)] $\aut(\mathcal M_1)\leq \aut(\mathcal M_2)=\aut(\mathcal M_3)
\leq \aut(\mathcal M_4)$, 
\item[(b)] $|\aut(\mathcal M_1)|=p^3$, 
\item[(c)] $|\aut(\mathcal M_2)|=2p^3$,
\item[(d)] $|\aut(\mathcal M_3)|=2p^3$, 
\item[(e)] $|\aut(\mathcal M_4)|=8p^3$.
\end{itemize}
\end{theorem}

\noindent{\bf Proof.} 
It is clear that the previously defined  permutations $t_{ab}$ and $\varphi$ 
are elements of each  automorphism group  $\aut(\mathcal M_i)$, hence $H$ 
is a subgroup of $\aut(\mathcal M_i)$ for each $1\leq i\leq 4$. 

\begin{itemize}
\item[(a)] The chain of inequalities     
$\aut(\mathcal M_1)\leq \aut(\mathcal M_2)\leq\aut(\mathcal M_3)\leq \aut(\mathcal M_4)$ 
follows directly from Theorem \ref{thm3}, simply by applying Galois correspondence 
to the lattice of coherent configurations and that of their corresponding   
automorphism groups. The equality $\aut(\mathcal M_2)=\aut(\mathcal M_3)$
will follow by inspection  of the group orders, to be accomplished in  
parts (c) and (d) below. 
\end{itemize}

Below we provide separate proofs for each of the claims  (b) through (e). 
We use the same methodology throughout. In each  proof,  $G$ will denote the group 
$\aut(\mathcal M_i)$, while $G_{[0,0],(0,0)}$ will denote  the stabilizer in $G$ of both the point
$[0,0]$ and the line $(0,0)$. The final result is obtained via manipulation of suitable 
elements of $G$ and application of the classical orbit-stabilizer lemma. 
Out of necessity, we will  introduce certain suitable permutations acting 
on $\mathcal P\cup\mathcal L$ that we have not already encountered.
  
\begin{itemize}
\item[(b)] We apply the orbit-stabilizer lemma to prove that $|\aut(\mathcal M_1)|\le p^3$.  
As we  already know that $H$  is a subgroup of $\aut(\mathcal M_1)$, the result will follow. 
(In fact, this will further show that 
$\aut(\mathcal M_1)\cong H \cong \mathbb Z_p^2\rtimes \mathbb Z_p$.) 

Denote $G:=\aut(\mathcal M_1)$ for brevity. First we claim that there is no 
automorphism that sends a point to a line. By way of contradiction, suppose 
$\alpha\in G$ sends the point $P_1$ to the line $r=(k,q)$,  $k,q\in \mathbb Z_p$.
Without loss of generality we may assume $P_1=[0,0]$, because $G$ acts transitively 
on $\mathcal P$. In such case, we must have $\mathcal P^{\alpha}=\mathcal L$ and 
$\mathcal L^{\alpha}=\mathcal P$, because of the relations $S_i$. 
Consider now the line $l=(0,0)$ and its image $l_1^{\alpha}=(u,v)$, $u,v\in\mathbb Z_p$. 
Clearly $(P_1,l_1)\in U_0$, whence $v=k\cdot u+q$. Now  consider the point $[1,0]$. 
Since $([0,0],[1,0])\in T_1$ and $([1,0],(0,0))\in U_0$, it follows that 
$[1,0]^{\alpha}=(k+1,q-u)$. Similarly, $(1,0)^{\alpha}=[u+1,q+k(u+1)]$. 
However $([1,0],(1,0))\in U_1$, and therefore $([1,0],(1,0))^{\alpha}\in U_1$.
But $([1,0],(1,0))^{\alpha}=((k+1,q-u),[u+1,q+k(u+1)])\in U_{-1}$, since 
$q+k(u+1)-(k+1)(u+1)-q+u=-1$, a contradiction for any odd prime $p$. 
This proves that $\mathcal P^G=\mathcal P$ and $\mathcal L^G=\mathcal L$, as claimed. 

Since $G$ is transitive on the points, $|[0,0]^G|=|\mathcal P|=p^2$. 
Let $G_{[0]}:=G_{[0,0]}$ be the stabilizer in $G$ of the point $[0,0]$. 
The points $[0,1],[0,2],\ldots,[0,p-1]$ are fixed by $G_{[0]}$, because 
they form unique pairs together with $[0,0]$ in the relations 
$S_1,S_2,\ldots,S_{p-1}$, respectively. As the line $(0,0)$ contains the point
$[0,0]$, there are at most  $p$ distinct images of $(0,0)$ under the action of $G_{[0]}$. 
However, it is easy to check that $(0,0)^{\varphi^i}=(-i,0)$ for $0\le i \le p-1$, 
which proves that $|(0,0)^{G_{[0]}}|=p$. 
 
Now let  $G_0:=G_{[0,0],(0,0)}$ be the stabilizer in $G$ of both $(0,0)$ and $[0,0]$. 
Then $G_0$ fixes all lines $(0,i)$ parallel to $(0,0)$, because $(0,i)$ forms a 
unique pair with $(0,0)$ in $S_i$. If we now consider an arbitrary point $[x,y]$ 
with $x\neq 0$, then its image under $G_0$ must be contained in  the line $(0,y)$.
Moreover, $([0,y],[x,y])\in T_x$ and for any $\pi\in G_0$ it follows that 
$([0,y],[x,y])^{\pi}=([0,y],[t,y])\in T_x$ for some  $t\in \mathbb Z_p$. 
This establishes that  $t=x$, and hence the point $[x,y]$ is fixed under $G_0$. 
Thus $G_0$ fixes all points and therefore all lines as well. By the orbit-stabilizer 
lemma, this yields $|G|=p^2\cdot p\cdot 1=p^3$ and so $G\cong H$ as desired. 

\item[(c)] By routine inspection, one can check that the mapping $\pi$ defined by  
$[x,y]\mapsto (x,-y-2x)$, $(x,y)\mapsto [x+2,-y]$ is an automorphism of $\mathcal M_2$. 
From this it follows that $G:=\aut(\mathcal M_2)$ is transitive on 
$V=\mathcal P\cup\mathcal L$, whence $|[0,0]^G|=2p^2$. As in the case above,  
we may again show that the line $(0,0)$ has $p$ distinct images under the action of 
the stabilizer $G_{[0]}$ of $[0,0]$. Let us again consider the stabilizer  $G_0$ of  
point $[0,0]$ and line $(0,0)$. For any point $[x,y]\in\mathcal P$ we have that 
$[x,y]^{G_0}\subseteq\{[x,y],[x,-y]\}$ because of the relations $S_i^*$ and $U_i^*$. 
If there exists some point $P_1=[x,y]$ ($x\neq 0$ and $y\neq 0$) for which
$|[x,y]^{G_0}|=2$, then we can write $y=k\cdot x$ for some $k\neq 0$, and 
$[x,kx]^h=[x,-kx]$, $[x,-kx]^h=[x,kx]$ which is equivalent to $(k,0)^h=(-k,0)$,
$(-k,0)^h=(k,0)$. But $((-k,0),(k,0))\in T_{2k}$, whence 
$((-k,0),(k,0))^h\in T_{2k}$ as well. This implies that $k=-k$, a contradiction 
since $k\neq 0$ and $p$~is odd. Hence $[x,y]$  is fixed by~$G_0$, and it follows that
all points and lines are fixed by~$G_0$. We conclude that  $|G|=2p^3$ as desired.

\item[(d)] One easily verifies that the  permutation $\pi$ of part (c) is also 
an automorphism of $\mathcal M_3$, i.e. $\pi\in G:=\aut(\mathcal M_3)$. 
Moreover, the initial steps of part (c) again establish that $|[0,0]^G|=2p^2$.
We consider once more the stabilizer  $G_0$ of  $(0,0)$ and $[0,0]$. Because of 
relations $T_i$, the points $[1,0], [2,0],\ldots, [p-1,0]$ are also stabilized 
by $G_0$. Let $\alpha\in G_0$. Since $T_i^{\alpha}=T_i$, we have 
$[x,y_1]^{\alpha}=[x,y_2]$ and $(k,q_1)^{\alpha}=(k,q_2)$, i.e. $\alpha$ does 
not change the first coordinate of a point or line. In particular, $\alpha$ 
preserves each parallel class of lines, hence it permutes the lines 
$(0,1),(0,2),\ldots,(0,p-1)$ amongst themselves. Note that the manner in which 
$\alpha$ permutes the lines $(0,i)$ uniquely determines the images of points 
under $\alpha$. Let $(1,0)^{\alpha}=(c,0)$ for some $c\in \mathbb Z_p^*$. 
Then for all $x,y\in\mathbb Z_p$, we necessarily have  $[x,y]^{\alpha}=[x,cy]$. 
But $((0,0),(1,0))\in T_1$, so $((0,0),(1,0))^{\alpha}\in T_1$, i.e. 
$((0,0),(c,0))\in T_1$, which forces $c=1$. Thus $G_0$ fixes all points as well 
as lines. It follows that $|G|=2p^3$ as desired. (In fact, we  confirmed  that   
$\aut(\mathcal M_2)\cong\aut(\mathcal M_3)\cong\langle t_{01},t_{10},\varphi,\pi\rangle$.) 

\item[(e)] Here we set $G:=\aut(\mathcal M_4)$, and consider  the permutations 
$\alpha,\beta$ defined by   $[x,y]^{\alpha}=[x,-y]$, $(k,q)^{\alpha}=(-k,-q)$ and 
$[x,y]^{\beta}=[-x,-y]$, $(k,q)^{\beta}=(k,-q)$.  
In a fashion similar to parts (b) and (c), one can verify that 
$G\cong \langle t_{10},t_{01},\varphi,\pi,\alpha,\beta\rangle$, whence   
$|G|=8p^3$ as claimed. 

\rightline{$\square$}

\end{itemize} 

As was mentioned earlier, the results presented in this section were obtained
as a theoretical generalisation of a tremendous number of computations, 
fulfilled with the aid of a computer for starting small values of prime numbers $p$.
In addition to the knowledge of combinatorial groups of automorphisms, we were obtaining
their ranks and were establishing that, as a rule, the proceeded association schemes
are non-Schurian. This,  together with Theorem~\ref{thm4}, allowed us to formulate 
our next result.

\begin{theorem}\label{thm5}
For $p>3$, $\mathcal M_1$, $\mathcal M_2$, $\mathcal M_3$, $\mathcal M_4$ 
are pairwise distinct non-Schurian association schemes.  
\end{theorem}

\medskip

\noindent{\bf Proof.} 
Recall our earlier derivation that the number of 2-orbits of 
$H=\mathbb Z_p^2\rtimes\mathbb Z_p$ on $V=\mathcal P\cup\mathcal L$ is $6p-2$.  
Thus  $\aut(\mathcal M_1)$  is of rank $6p-2$, while the rank of $\mathcal M_1$ is $3p-1$.  
This  proves  that~$\mathcal M_1$ is non-Schurian for $p\ge 3$.  

For the association schemes $\mathcal M_2$ and $\mathcal M_3$, we consider the permutation 
$\pi\in \aut(\mathcal M_2)=\aut(\mathcal M_3)$ introduced in part (c) of Theorem~\ref{thm4}.  
As the result of $\pi$, we obtain the following 2-orbits:  
$A_i\cup C_{-i}$, $B_j\cup D_j$ and $E_i\cup F_{-i}$, for $i\in\mathbb Z_p$ and 
$j\in \mathbb Z_p^*$. This proves that $\aut(\mathcal M_2)=\aut(\mathcal M_3)$ 
is of rank $3p-1$. As the ranks of  $\mathcal M_2$ and~$\mathcal M_3$ are $2p$ and 
$p+3$ respectively, we conclude that $\mathcal M_2$ and~$\mathcal M_3$ are 
non-Schurian for $p\ge 3$. 
However, one can check that $\mathcal M_2$ and~$\mathcal M_3$ coincide when $p=3$.
 
Finally, as a result of the permutations $\alpha,\beta$ introduced in part (e)
of Theorem \ref{thm4}, it is easy to see that the 2-orbits of $\aut(\mathcal M_4)$
are $A_i\cup A_{-i}\cup C_i\cup C_{-i}$, $B_j\cup B_{-j}\cup D_j\cup D_{-j}$ and $
E_i\cup E_{-i}\cup F_i\cup F_{-i}$, for $i\in \{0,1,\ldots,\frac{p-1}{2}\}$ and 
$j\in\{1,2,\ldots,\frac{p-1}{2}\}$. Thus the rank of $\aut(\mathcal M_4)$ is equal 
to $\frac{3p+1}2$.  As the rank of  $\mathcal M_4$ is $\frac{p+7}2$, we conclude 
that $\mathcal M_4$ is non-Schurian for $p>3$.
\rightline{$\square$}

\medskip

\noindent{\bf Remark 4. } Note that when $p=3$, we get only two non-Schurian 
association schemes, namely $\mathcal M_1$ and $\mathcal M_2$. Indeed, when  
$p=3$ the association schemes $\mathcal M_2$ and $\mathcal M_3$ coincide,   
while $\mathcal M_4$ is  a Schurian association scheme  
(here $\frac{3\cdot 3+1}2=\frac{3+7}2=5$). 
These two association schemes are of order $2\cdot 3^2=18$, and according to
the catalogue of small association schemes of Hanaki and Miyamoto~\cite{hm}, 
we see that our constructions cover all non-Schurian association
schemes of order~$18$.

\section{A second model of $\mathcal M$}

According to our intentionally chosen genetic style of presentation, we now consider 
our favourite objects once again, in a sense, from scratch.

Recall that in the first part of the text our starting object was the biaffine plane 
$\mathcal B_p$, and our preferred descriptive language was geometric in nature.  
Out of necessity, suitable permutations appeared in an ad hoc manner. 

In this section we describe a new second model of $\mathcal M$.  (Isomorphism with the first 
model will be shown later on.) The main advantage of this second model is its purely 
algebraic flavour. As a result, many of our previous claims will now get more transparent
proofs, probably more preferable  for readers with developed algebraic tastes.

\subsection{Initial definitions}

We start with description of the second model of $\mathcal M$. 

Let
\[
V_1=\mathbb Z_p^2=\{ (1, x_1, x_2) \mid x_1, x_2\in \mathbb Z_p\},
\]
\[
V_2=(\mathbb Z_p^2)^{dual}=\left\{ \left( \begin{array}{c} x_1 \\ x_2 \\ -1\end{array}\right)
\left| \right. x_1, x_2\in \mathbb Z_p\right\}.
\]
Define a natural scalar product between $V_1$ and $V_2$:
\[
(1, x_1, x_2)\left( \begin{array}{c} y_1 \\ y_2 \\ -1\end{array}\right)=y_1+x_1y_2-x_2.
\]
Let 
\[
H'=\left\{g_{abc}=\left( \begin{array}{ccc} 1 & a & b+ac \\ 0 & 1 & c \\ 0 & 0 & 1 \end{array}\right)
\left|\right. a,b,c \in \mathbb Z_p\right\}.
\]
Matrix $g_{abc}$ is invertible, and 
\[
g^{-1}_{abc}=\left( \begin{array}{ccc} 1 & -a & -b \\ 0 & 1 & -c \\ 0 & 0 & 1 \end{array}\right).
\]
Clearly, the set $H'$ together with the operation of matrix-multiplication forms a group. 
The multiplication in $H'$ is given by $g_{abu}\cdot g_{cdv}=g_{a+c,b+d-uc,u+v}$,
and this group is well known under the name \emph{Heisenberg group modulo $p$,}
see e.g. \cite{bw}. 

Define an action of $H'$ on $\Omega=V_1\cup V_2$ by: 
\(
x^g=\begin{cases} x\cdot g & \text{if } x\in V_1 \\
g^{-1}\cdot x & \text{if } x \in V_2, \end{cases}
\, \text{for all } g\in H'. 
\)

If we take arbitrary $x_1,x_2,y_1,y_2\in\mathbb Z_p$, then 
for $a=y_1-x_1$, $b=y_2-x_2$, $c=0$ matrix $g_{abc}$ sends 
$(1,x_1,x_2)\in V_1$ to $(1,y_1,y_2)\in V_1$, 
while $g_{cab}$ sends $(x_1,x_2,-1)^T\in V_2$ to
$(y_1,y_2,-1)^T\in V_2$, therefore $H'$ acts transitively on $V_1$
and $V_2$ as well. 

The action of $H'$ preserves the scalar product: $x^g\cdot y^g =
(x\cdot g)\cdot(g^{-1}\cdot y)=x\cdot gg^{-1}\cdot y=x\cdot y$. 

\begin{proposition}\label{prop6}
Groups $H$ and $H'$ are isomorphic.
\end{proposition}

\noindent{\bf Proof.} We claim that $\Phi:H\to H'$, 
$h_{a,b,u}\mapsto g_{a,b+au,-u}$ is a group isomorphism.  First, we have
\begin{align*}
\Phi(h_{a,b,u}\circ h_{c,d,v})&=\Phi(h_{a+c,b+d-av,u+v})=
g_{a+c,b+d-av+(a+c)(u+v),-(u+v)}=\\
&=g_{a,b+au,-u}\cdot g_{c,d+cv,-v}=
\Phi(h_{a,b,u})\cdot\Phi(h_{c,d,v}).
\end{align*}
Moreover, $\Phi$ is invertible: $\Phi^{-1}(g_{\alpha,\beta,\gamma})=
h_{\alpha,\beta+\alpha\gamma,-\gamma}$. Thus,  $\Phi$ is a group
isomorphism from $H$ to $H'$ as claimed.

\rightline{$\square$}

Proposition \ref{prop6} allows us to identify groups $H$ and $H'$, which is 
why we shall henceforth ascribe the notation $H$ to both groups. 
However, the reader is advised to keep  in mind the group $H'$  
as it appears in this section.

\subsection{Orbits of $H$ on $\Omega^2$}

The isomorphism between groups $H$ and $H'$  established in  
Proposition \ref{prop6} allows an alternate description of the 2-orbits of the 
biaffine coherent configuration $\mathcal M$, and correspondingly, an alternate proof  
of Proposition \ref{prop1}.

Not attempting to put such a proof into the framework of formal propositions, 
we nevertheless outline   the necessary arguments and final
formulations for the reader's benefit. 

\medskip

\noindent{\bf Orbits on $V_1\times V_1$:}

Let $g=g_{abc}\in H$ and $P_1=(1,x_1,x_2)$,
$P_2=(1,y_1,y_2)$, $P_3=(1,u_1,u_2)$, $P_4=(1,v_1,v_2)\in V_1$. 
Then $(P_1,P_2)^g=(P_3,P_4)$ if and only if
\begin{eqnarray*}{}
u_1 &=& x_1+a \\ v_1&=&y_1+a \\ u_2 &=& c\cdot x_1+x_2+b+a\cdot c \\
v_2 &=& c\cdot y_1+y_2+b+a\cdot c.
\end{eqnarray*}

We can see that if $(P_1,P_2)$ and $(P_3,P_4)$ belong to the same 
orbit, then necessarily $y_1-x_1=v_1-u_1$, $a=u_1-x_1=v_1-y_1$, and
$v_2-u_2=c(y_1-x_1)+(y_2-x_2)$. 
\begin{itemize}
\item[(i)] If $y_1-x_1=0$, then it is necessary to have $v_2-u_2=y_2-x_2$, 
and by choosing $c=0$ and $b=u_2-x_2$ we are getting a suitable $g$ sending
$(P_1,P_2)$ to $(P_3,P_4)$.
\item[(ii)] If $y_1-x_1=k\in\mathbb Z_p^*$, then 
we can choose $c=k^{-1}(v_2-u_2+x_2-y_2)$, $b=u_2-cx_1-x_2-ac$
and we get $(P_1,P_2)^g=(P_3,P_4)$.
\end{itemize}

Hence we have got two types of orbits, say $A_k$ and $B_k$ on $V_1\times V_1$:
\begin{itemize}
\item $(P_1,P_2)\in A_k$ if and only if $y_1=x_1$ and $y_2-x_2=k$, where $k\in
\mathbb Z_p$.
\item $(P_1,P_2)\in B_k$ if and only if $y_1-x_1=k$, where $k\in
\mathbb Z_p^*$.
\end{itemize} 

\noindent{\bf Orbits on $V_2\times V_2$:}

Let $g=g_{abc}\in H$ and $l_1=(x_1,x_2,-1)^T$,
$l_2=(y_1,y_2,-1)^T$, $l_3=(u_1,u_2,-1)^T$, $l_4=(v_1,v_2,-1)^T\in V_2$. 
Then $(l_1,l_2)^g=(l_3,l_4)$ if and only if
\begin{eqnarray*}{}
u_1 &=& x_1-ax_2+b \\ v_1&=&y_1-ay_2 \\ u_2 &=& x_2+ c \\
v_2 &=& y_2+ c.
\end{eqnarray*}

We can see that if $(l_1,l_2)$ and $(l_3,l_4)$ belong to the same 
orbit, then necessarily $y_2-x_2=v_2-u_2$, $c=u_2-x_2=v_2-y_2$, and
$v_1-u_1=a(x_2-y_2)+(y_1-x_1)$. 
\begin{itemize}
\item[(i)] If $y_2-x_2=0$, then it is necessary to have $v_1-u_1=y_1-x_1$, 
and by choosing $a=0$ and $b=v_1-y_1$ we obtain a suitable $g$ sending
$(l_1,l_2)$ to $(l_3,l_4)$.
\item[(ii)] If $x_2-y_2=k\in\mathbb Z_p^*$, then 
we can choose $a=k^{-1}(v_1-u_1+x_1-y_1)$, $b=v_1-y_1+ay_2$
and we get $(l_1,l_2)^g=(l_3,l_4)$.
\end{itemize}

Hence we get two types of orbits, say $C_k$ and $D_k$, on $V_2\times V_2$:
\begin{itemize}
\item $(l_1,l_2)\in C_k$ if and only if $y_2=x_2$ and $y_1-x_1=k$, where $k\in
\mathbb Z_p$.
\item $(l_1,l_2)\in D_k$ if and only if $y_2-x_2=k$, where $k\in
\mathbb Z_p^*$.
\end{itemize} 

\noindent{\bf Orbits on $V_1\times V_2$:}

Let $P_1=(1,x_1,x_2)$, $P_2=(1,y_1,y_2)\in V_1$ and
$l_1=(u_1,u_2,-1)^T$, $l_2=(v_1,v_2,-1)^T\in V_2$.

Since our scalar product is preserved by $H$, it follows that $(P_1,l_1)$ and 
$(P_2,l_2)$ belong to different orbits when $P_1\cdot l_1\neq P_2\cdot l_2$. 
Now suppose $P_1\cdot l_1=P_2\cdot l_2$, i.e. 
\begin{equation}
y_2=v_1+y_1v_2+x_2-u_1-x_1u_2. 
\end{equation}
We will show that in this case there exists 
$g\in H$ such that $(P_1,l_1)^g=(P_2,l_2)$. 
Necessarily, 
\begin{eqnarray*}
y_1 &=& x_1+a \\
y_2 &=& x_2+cx_1+b+ac \\
v_1 &=& u_1-au_2+b \\
v_2 &=& u_2+c.
\end{eqnarray*}
Choosing $a=y_1-x_1$, $c=v_2-u_2$, $b=v_1-u_1+au_2$, all equations are fulfilled 
and therefore $P_1^g=P_2$, $l_1^g=l_2$. Thus the orbits $E_k$ of $H$ on $V_1\times V_2$
are formed by pairs with the same value of scalar product: $(P,l)\in E_k \iff P\cdot l=k$,
where $k\in\mathbb Z_p$.

\medskip

\noindent{\bf Orbits on $V_2\times V_1$:}
 
Similar to the previous case, we obtain orbits $F_k$: 
 $(l,P)\in F_k \iff P\cdot l = -k$. 

\medskip

Altogether, we get the same six types of orbits of $H$ on $\Omega\times\Omega$ 
with which we are already familiar:
\begin{eqnarray*}
(P_1,P_2)\in A_k &\iff& y_1=x_1 \text{ and } y_2-x_2=k\in\mathbb Z_p \\
(P_1,P_2)\in B_k &\iff& y_1-x_1=k\in\mathbb Z_p^* \\ 
(l_1,l_2)\in C_k &\iff& u_2=v_2 \text{ and } v_1-u_1=k\in\mathbb Z_p \\
(l_1,l_2)\in D_k &\iff& v_2-u_2=k\in\mathbb Z_p^* \\ 
(P_1,l_1)\in E_k &\iff& P_1\cdot l_1=k\in\mathbb Z_p \\
(l_1,P_1)\in F_k &\iff& P_1\cdot l_1=-k\in\mathbb Z_p, 
\end{eqnarray*}
where $P_1=(1,x_1,x_2)$, $P_2=(1,y_1,y_2)\in V_1$ and
$l_1=(u_1,u_2,-1)^T$, $l_2=(v_1,v_2,-1)^T\in V_2$.

\section{Toward a theoretical understanding of the algebraic groups for the
detected coherent configurations}

At this stage we are ready to discuss how, in principle, one may describe the
full algebraic group of a prescribed coherent configuration $\mathcal W$, in
particular, of an association scheme. 

Recall that $\aaut(\mathcal W)$ acts faithfully on the set of relations from $\mathcal W$, 
preserving the tensor $\mathcal T$ of structure constants of $\mathcal W$. 

Roughly speaking, our methodology may be described as follows:
\begin{itemize}
\item detect a few permutations, say $\varphi_1,\varphi_2,\ldots,\varphi_k$, acting on
the set of relations of $\mathcal W$;
\item check that each permutation $\varphi_i$, $1\leq i\leq k$, preserves $\mathcal T$;
\item describe group $K=\langle \varphi_1,\varphi_2,\ldots,\varphi_k\rangle$  as abstract
group, and establish its orbits on the set of relations of  $\mathcal R$;
\item prove, using ad hoc arguments, that each element $\varphi\in\aaut(\mathcal W)$ may be
expressed as a suitable concatenation $x_1x_2\ldots x_l$, where  
$x_j\in\{\varphi_1,\varphi_2,\ldots,\varphi_k\}$ for $1\leq j\leq l$. 
\end{itemize}
Note that in general, the final step seems quite  sophisticated. 
Nevertheless, for small numbers~$k$ of generators, the tricks we employed appear
to be fairly visible 
and reasonably natural. 

There is a definite sense to stress that the problem of calculating $\aaut(\mathcal W)$ is
a relatively new kind of activity in AGT. The paper \cite{km} was among the first in which 
such types of reasonings were considered. In a number of other 
publications the calculation was achieved with the aid of a computer. One of the foremost goals of the  
current text is to introduce some helpful tricks and theoretical reasonings which allow one,
in principle, to attack the problem with a sufficient level of rigour.  
Thus, we will first demonstrate the suggested methodology on the association scheme $\mathcal M_1$. 

\begin{proposition}\label{prop7}
Permutations
\begin{eqnarray*}
\alpha &=& (T_1,T_{\omega},T_{\omega^2},\ldots,T_{\omega^{p-2}}), \\
\beta &=& (S_1,S_{\omega},S_{\omega^2},\ldots,S_{\omega^{p-2}})
(U_1,U_{\omega},U_{\omega^2},\ldots,U_{\omega^{p-2}}),
\end{eqnarray*}
acting on the set of relations of the association scheme $\mathcal M_1$, preserve its tensor
$\mathcal T_1$ of structure constants. 
Here $\omega$ is a primitive element of $\mathbb Z_p$ regarded as a field 
(in other words, $\omega$ is a generator of the multiplicative group $\mathbb Z_p^*$). 
\end{proposition}

\noindent{\bf Proof.}
We refer to the Appendix 1, where the tensor $\mathcal T_1$ is presented with the aid
of four tables. 
Permutation $\alpha$ acts only on the set of relations of type $T_i$ as
$T_i^{\alpha}=T_{i\omega}$, therefore
it is enough to check what is happening with the intersection numbers
related to this type of relations. For example,
\[
c_{ti,tj}^0=p\cdot \delta_{i,j'}, \quad 
c_{ti^{\alpha},tj^{\alpha}}^0=c_{ti\omega, tj\omega}^0=p\cdot
\delta_{i\omega,(j\omega)'}. 
\]
So we need to show that $\delta_{i,j'}=\delta_{i\omega,(j\omega)'}$. 
Our relations satisfy $i'=-i$, hence the equation $i=-j$ is equivalent 
to $i\omega=-j\omega$, in which case we are done. 
In the case of intersection numbers of the type $c_{ti,tj}^{tk}$ we need to
show that $\delta_{i+j,k}=\delta_{i\omega+j\omega,k\omega}$, which is obvious.
Using similar reasonings, one proves   
that
$\beta$ preserves $\mathcal T_1$ as well. 

\rightline{$\square$}

We now observe that both $\alpha$ and $\beta$ are cyclic permutations of
order $p-1$. Also, we easily check that $\alpha$ and $\beta$ commute by making 
a straightforward comparison of  $\alpha\beta$ with $\beta\alpha$. This implies that
$\langle\alpha,\beta\rangle\cong\mathbb Z_{p-1}^2$. Denote by $K$ this group 
$\mathbb Z_{p-1}^2$ in its action on the relations of  $\mathcal M_1$. 

\begin{corollary}\label{cor8}
$\aaut(\mathcal M_1)\geq K$.
\end{corollary}

Our next step is  simply to prove that each $\varphi\in\aaut(\mathcal M_1)$ 
indeed belongs to $K$. 

\begin{theorem}\label{thm9}
$\aaut(\mathcal M_1)\cong \mathbb Z_{p-1}^2$. 
\end{theorem}

\noindent{\bf Proof.} 
Let $\varphi\in\aaut(\mathcal M_1)$. 
The only relation $i$ in $\mathcal M_1$ for which $c_{ii}^i=1$ is $i=R_0$, 
therefore $R_0^{\varphi}=R_0$. Since $c_{si,si'}^0=1$ and there are no
other 1's among the values of type $c^0$, we have that $S_i^{\varphi}=S_j$ for some $j$. 
Moreover, $S_{i'}^{\varphi}=S_{j'}$. Just using values of type $c^0$, similar considerations 
lead to the following observations:  for all $i$ there exists $j$ such that 
$T_i^{\varphi}=T_j$, and for all $l$ there exists $m$ such that $U_l^{\varphi}=U_m$. 
In particular, $U_0^{\varphi}=U_0$.

Now let $S_1^{\varphi}=S_{\omega}$, where $\omega$ is a primitive element of $\mathbb Z_p^*$.
Then by induction on $i$,  $c_{s1,si}^{si+1}=1$ implies $S_{i+1}^{\varphi}=S_{\omega(i+1)}$. 
Similarly $T_i^{\varphi}=T_{\omega}$ implies $T_i^{\varphi}=T_{\omega\cdot i}$.
However $S_i^{\varphi}=S_{\omega\cdot i}$, and $c_{u0,ui}^{si}=p=c_{u0,uj}^{s\omega i}$,
therefore $j=\omega i$, i.e., $U_i^{\varphi}=U_{\omega i}$.

Letting $\omega$ be a primitive element of $\mathbb Z_p$, the above computations 
establish that $\aaut(\mathcal M_1)\leq \langle\alpha,\beta\rangle$.  
Together with Corollary \ref{cor8}, we conclude have $\aaut(\mathcal M_1)=
\langle \alpha,\beta\rangle$ from which follows  $\aaut(\mathcal M_1)\cong
\mathbb Z_{p-1}^2$. 

\rightline{$\square$}

\section{Algebraic group of the coherent configuration $\mathcal M$}

In the previous section we gave a description of the group  $\aaut(\mathcal M_1)$. 
We hope that our presentation  fulfilled its mission: to create for the
reader a first acquaintance with this significant concept. 

In this section we will describe the group $\aaut(\mathcal M)$. 
Knowledge of the group  will play a critical role in achieving our further goals. 
Once more we wish to stress significant features of the used technology. Indeed, 
at the first stage the description of the order of $\aaut(\mathcal M)$ was 
obtained for initial values of the prime parameter $p$. Again aided by {\sf GAP}, 
we next made educated guesses at the structure of these groups for the obtained 
orders.  Finally, an evident formulation of a theorem became transparent to us. 
From this we created a pedestrian proof, a brief outline of which is given below.

Thus let us start with the following permutations on the set of relations of the coherent
configuration $\mathcal M$.  (We assume $\omega$ is a primitive element of $\mathbb Z_p^*$.)
\begin{eqnarray*}
g_1 &=& (A_0,C_0)(E_0,F_0)\prod_{i=1}^{p-1}(A_i,C_i)(B_i,D_i)(E_i,F_i), \\
g_2 &=& (E_0,E_{p-1},E_{p-2},\ldots,E_2,E_1)(F_0,F_1,F_2,\ldots,F_{p-1}), \\
g_3 &=& (A_1,A_{\omega},A_{\omega^2},\ldots,A_{\omega^{p-2}})(C_1,C_{\omega},
\ldots,C_{\omega^{p-2}})(E_1,E_{\omega},\ldots,E_{\omega^{p-2}})(F_1,F_{\omega},
\ldots,F_{\omega^{p-2}}), \\
g_4 &=& (B_1,B_{\omega},B_{\omega^2},\ldots,B_{\omega^{p-2}}), \\
g_5 &=& (D_1,D_{\omega},D_{\omega^2},\ldots,D_{\omega^{p-2}}). 
\end{eqnarray*}

Permutation $g_1$ is an involution which corresponds to the duality between points
and lines in $\mathcal B_p$ (it interchanges the roles of points and lines). 
Permutation $g_2$ is of order $p$, while the remaining permutations are of order $p-1$. 

\begin{theorem}\label{thm10}
The group $\aaut(\mathcal M)$ is of order $2p(p-1)^3$ and  
$$ \aaut(\mathcal M)\cong\langle g_1,g_2,g_3,g_4,g_5\rangle
\cong (\mathbb Z_{p-1}^2 \rtimes \mathbb Z_2)\times \AGL(1,p).$$
\end{theorem}

\noindent{\bf Proof.} 
First, we have to check that each of the permutations $g_1,\ldots,g_5$ preserves
the tensor of structure constants. This job can be done by hand, although
we admit it  is a bit tedious. Nevertheless, we suggest that inspection of at 
least one of the presented permutations may serve as a helpful exercise for 
the reader. Thus, we have that $\langle g_1,g_2,g_3,g_4,g_5\rangle\leq \aaut(\mathcal M)$.
In order to see that $\aaut(\mathcal M)\cong\langle g_1,g_2,g_3,g_4,g_5\rangle$ it is
sufficient to compute and compare the orders of  $\langle g_1,g_2,g_3,g_4,g_5\rangle$ 
and $\aaut(\mathcal M)$. The order of $\aaut(\mathcal M)$ may be determined by 
repeated application of the orbit-stabilizer lemma. The order of 
$\langle g_1,g_2,g_3,g_4,g_5\rangle$ will be derived by us presently. 

Clearly $\{g_2,g_3\}$ is a standard set of generators for the affine linear group 
over the finite field $\mathbb F_p$ of order $p$, that is, 
$\langle g_2,g_3\rangle\cong \AGL(1,p)$. Also it is easy to see that 
$\langle g_4,g_5\rangle\cong \mathbb Z_{p-1}\times\mathbb Z_{p-1}=\mathbb Z_{p-1}^2$ and 
$\langle g_1,g_4,g_5\rangle\cong \mathbb Z_{p-1}^2\rtimes_{\phi}\mathbb Z_2$,
where $\phi$ interchanges $g_4$ with $g_5$. By routine computation one can  
further show that $\langle g_1,g_4,g_5\rangle\unlhd \langle g_1,g_2,g_3,g_4,g_5\rangle$, 
$\langle g_2,g_3\rangle\unlhd \langle g_1,g_2,g_3,g_4,g_5\rangle$, and
$\langle g_1,g_4,g_5\rangle\cap\langle g_2,g_3\rangle=\{e\}$. 
Thus $\langle g_1,g_4,g_5\rangle\times\langle g_2,g_3\rangle\cong
\langle g_1,g_2,g_3,g_4,g_5\rangle$, which confirms that 
$\langle g_1,g_2,g_3,g_4,g_5\rangle\cong(\mathbb Z_{p-1}^2\rtimes\mathbb Z_2)\times \AGL(1,p)$.
Finally, $|\langle g_1,g_2,g_3,g_4,g_5\rangle|=2p(p-1)^3$.

\rightline{$\square$}

\section{Detected non-Schurian association schemes as algebraic mer\-gings}

As we have seen in the previous section, the algebraic group $\aaut(\mathcal M)$
is of order $2p(p-1)^3$. It is easy to check that this group has four orbits of 
length $2$, $2p-2$, $2p-2$ and $2p$, respectively, on the set of relations 
of $\mathcal M$. These orbits are:
\[
\{A_0,C_0\}, \{A_1,\ldots,A_{p-1},C_1,\ldots,C_{p-1}\}, \{B_1,\ldots,B_{p-1},
D_1,\ldots,D_{p-1}\},\]
\[
\{E_0,E_1,\ldots,E_{p-1},F_0,F_1,\ldots,F_{p-1}\}.
\]
Knowledge of the algebraic group of automorphisms of a coherent configuration 
is important for constructing algebraic mergings. It is known (see \cite{km}) 
that if we take the orbits of a subgroup~$H$ of $\aaut(\mathcal W)$, then the 
algebraic merging with respect to the partition of the relations into orbits 
of $H$, leads to a coherent configuration. In particular, if $\mathcal W$ is 
an association scheme, then its algebraic merging is also an association scheme. 

\subsection{New interpretation of the main non-Schurian association schemes}

Here we present one more result, which was initially obtained via plausible
reasonings based on observation of computer aided data  for diverse values of $p$.
Its essential feature is that the proof becomes almost trivial, provided our 
suggested methodology of the use of algebraic groups is exploited in the correct manner.

\begin{proposition}\label{prop11}
The association schemes $\mathcal M_1, \mathcal M_2, \mathcal M_3,$ and
$\mathcal M_4$ appear as algebraic mergings of $\mathcal M$. 
\end{proposition}

\noindent{\bf Proof.} 
Let  $q=(p-1)/2$.  Then the algebraic mergings corresponding to the subgroups
$K_1=\langle g_1\rangle$, $K_2=\langle g_1, g_3^q\rangle$,
$K_3=\langle g_1, g_3\rangle$ and $K_4=\langle g_1, g_3, g_4^q, g_5^q\rangle$
of $\aaut(\mathcal M)$ lead to the association schemes $\mathcal M_1$, $\mathcal M_2$,
$\mathcal M_3$ and $\mathcal M_4$, respectively. 

\rightline{$\square$}

This proposition provides an alternate proof that the color graphs
defined in Section~6 correspond to association schemes. 
Recall that their ranks are $3p-1$, $2p$, $p+3$ and $(p+7)/2$, respectively. 

We are already aware that for $p>3$  these algebraic mergings are non-Schurian.

\subsection{Experimental results with color and algebraic automorphisms}

Computer experiments performed over a sufficiently large interval of initial values of the prime parameter $p$
show that $|\caut(\mathcal M)|=2p^4(p-1)^2$. 
This implies that the order of the quotient group $\caut(\mathcal M)/\aut(\mathcal M)$
is $\frac{2p^4(p-1)^2}{p^3}=2p(p-1)^2$. Comparing this   with the order of 
the group $\aaut(\mathcal M)$, we see that the index of the quotient group in 
$\aaut(\mathcal M)$ is $p-1$. Thus for all $p\geq 3$ there exist proper algebraic mergings.

It was at this particular stage of our project  that we became greatly enthused.  
Indeed, we fully realized  that the existence of proper algebraic mergings  
would create for us  a most attractive research agenda that could potentially 
guide us  to promising new discoveries.

In particular, for a few initial values of $p$ we arranged a full search of 
subgroups of $\aaut$, constructions of corresponding mergings, and an 
investigation of the properties of the resulting coherent configurations and 
association schemes. 

At this stage we are not yet prepared to transform our observations into 
the rigourous platform of proved mathematical claims. Indeed, our intentions 
are to do this in forthcoming publications (see discussion at the end of the text). 
Nevertheless, a portion of them will be considered in more detail in Section~14.

\section{Links to known combinatorial structures}

In this section we mention a few well known combinatorial structures, especially 
graphs, which are somehow related to the color graphs we have explored in this paper.
 
The graph defined by color $U_0$ for $p=3$ goes by  the name \emph{Pappus graph}. 
It is the incidence graph of the \emph{Pappus configuration} (see e.g. \cite{cx}). 

The so-called \emph{McKay-Miller-\v Sir\'a\v n graphs} $H_p$ are well known
in the degree/diameter problem (see \cite{ms} for a survey) which falls within 
the realm of extremal graph theory. They were defined originally in \cite{mm}, 
later on \v Siagiov\'a \cite{si} and finally Hafner \cite{ha} gave simplified 
constructions. We will define them according to Hafner's description. Though he 
described these graphs in terms of biaffine planes for arbitrary prime powers,  
for our purposes it suffices to restrict our attention to the case of odd primes. 
So let $p$ be an odd prime and put 
$V_p=\mathbb Z_2\times\mathbb Z_p\times\mathbb Z_p$ as the vertex set of $H_p$. 
Let $\omega$ be a primitive  element of  $\mathbb Z_p^*$. 
Since $p$ is an odd prime, there are two possibilities:
\begin{itemize}
\item If $p=4r+1$ then define $X=\{1,\omega^2,\omega^4,\ldots,
\omega^{p-3}\}$, $X'=\{\omega,\omega^3,\ldots,\omega^{p-2}\}$.
\item If $p=4r+3$ then define $X=\{\pm 1,\pm\omega^2,\ldots,
\pm\omega^{2r}\}$, $X'=\{\pm\omega,\pm\omega^3,\ldots,
\pm\omega^{2r+1}\}$.
\end{itemize}

Adjacency in  $H_p$ is defined as follows:
\begin{eqnarray*}
(0,x,y)  \text{ is adjacent to } (0,x,y') & \text{ if and only if } & y-y'\in
X, \\
(1,k,q) \text{ is adjacent to } (1,k,q') & \text{ if and only if } & q-q'\in
X', \\
(0,x,y) \text{ is adjacent to } (1,k,q) & \text{ if and only if } & y=kx+q.
\end{eqnarray*}

From this description it is clear that McKay-Miller-\v Sir\'a\v n graphs may 
be obtained as a suitable merging of relations of $\mathcal M$. Specifically:
$$ H_p = E_0\cup F_0 \cup \bigcup_{i\in X}A_i \cup\bigcup_{j\in X'}C_j. $$
In particular, $H_5$ is the well known Hoffman-Singleton graph
\cite{hs}.

The McKay-Miller-\v Sir\'a\v n graphs are currently the best known solutions 
to the degree/diameter problem for diameter 2 and valency $(3q-1)/2$, where 
$q$ is a prime power. They still play a significant role in newer 
constructions that give denser graphs, i.e. graphs that are closer to the Moore bound.

Recently, motivated by the success of McKay-Miller-\v Sir\'a\v n graphs,
authors in \cite{bm} investigated how to extend the graph corresponding to
$U_0$ in order to find better constructions in the degree/diameter problem. 
We plan to consider this operation of extension with more details in the future, 
because it fits quite well the language of relations in $\mathcal M$  
developed in our paper.

Another family of graphs which may be defined in our terminology is the 
family of \emph{Wenger graphs} $W_1(p)$ introduced in \cite{we}. They have 
$2p^2$ vertices with edge set coinciding with the relation $U_0=E_0\cup F_0$. 
In other words, they are flag graphs of the biaffine plane. The Wenger graphs 
were studied for their extremal properties and belong to a richer family 
of graphs defined by systems of equations on the coordinates of points and lines. 
For more details, see \cite{lu,lw,vi,wo}.

\section{Some association schemes of small rank, appearing as mergings of $\mathcal M$}
 
Recall that the main result of this paper is the discovery and investigation 
of four infinite families of non-Schurian association schemes, which arise as 
mergings of the coherent configuration $\mathcal M$ defined on $2p^2$ elements 
of the classical biaffine plane of odd prime order $p$. For each of these four 
families, rank of the schemes grows with increasing $p$.

From the earliest inception of AGT, special attention has been paid to association
schemes of small rank. The smallest possible (non-trivial) rank is 3, which 
corresponds to strongly regular graphs in the symmetric case, and doubly regular 
tournaments in the non-symmetric case. Applying special efforts, one can attempt 
to prove that in our case, for $p>5$, such primitive objects cannot appear.
Nevertheless, computer experiments created for us evidence that already for small values
of $p$ we are getting non-Schurian mergings of constant low ranks 5 and 6. This immediately
created a new challenge for us: how to justify existence of such schemes for arbitrary
values of $p$. 

The corresponding results are more fresh, with a portion of them having been announced 
by the author M.K.\ at the conference ``Modern trends in AGT'' held at Villanova 
University in June 2014. (See Section 15.9 for a link to the slides of this announcement.) 

In this section we provide an outline of the report of  these results. It is given in the  
form of research announcement, that is, we are not aiming  to give a full formulation 
or justification of our claims.
 
\subsection{Starting rank 6 scheme}

We begin with a brief discussion of our technology, which we feel may hold independent
interest due to its innovative combination of  computation  and further reasonings. 

Recall that for a few small values of $p$ a full enumeration of all coherent subalgebras 
was obtained (see Section 15). Further analysis of these  results suggested to us 
that it might be possible to obtain a non-Schurian merging association scheme of 
rank 6 for every value of $p$. Moreover, we found evidence to support that such a 
scheme might appear as a suitable algebraic merging. This led us to consider the 
subgroup $K$ of $\aaut(\mathcal M)$, where  $K=\langle g_1,g_3,g_4^2,g_5^2\rangle$. 
(Here we are following the notation introduced in Section 11.) 

Our next step was to apply our algebraic group $K$ to $\mathcal M$ for all odd primes 
$p\leq 19$. In each case we  constructed the corresponding algebraically merged 
association scheme and investigated its main properties.
For the reader's convenience a summary of our results is presented in Appendix 2.

A more careful theoretical analysis of our observations allowed us to reach the general picture
for arbitrary primes $p$. 

\medskip

\noindent{\bf Announcement 1.} For all odd primes $p$ there exists a non-Schurian rank 6 
algebraic merging~$\mathcal N_6$ of the master coherent configuration $\mathcal M$.
The group $\aut(\mathcal N_6)$ is a transitive rank 8 group of order $\frac12(p-1)^2p^3$.
The group $\caut(\mathcal N_6)$ has twice larger order. The  group $\aaut(\mathcal N_6)$
has order~2 and thus coincides with the group $\caut(\mathcal N_6)/\aut(\mathcal N_6)$.
The scheme $\mathcal N_6$ is commutative. It is non-symmetric when $p\equiv 3\pmod 4$
and symmetric when $p\equiv 1\pmod 4$. 

At the next stage we were able to describe the full tensor of structure constants 
of the scheme~$\mathcal N_6$. The results depend on the modulo 4 congruence class  
of~$p$, hence they are presented as two separate cases in the slides 
of M.K.\ mentioned above.   
 
Moreover, based on the use of classical techniques in AGT and relying on the knowledge of structure
constants, we were able to describe the characteristic polynomials and spectra of the basic
graphs of $\mathcal N_6$.

\medskip

\noindent{\bf Announcement 2.} Denote by $\Lambda_i$ the spectrum of the basic graph 
$\mathcal N_{6,i}$ in $\mathcal N_6$. Then:
\begin{eqnarray*}
\Lambda_1 &=& \{ 1 ^{2p^2} \}\\
\Lambda_2 &=& \{  p-1^{2p}, -1^{2p^2-2p}\} \\
\Lambda_3 = \Lambda_4 &=& \begin{cases}
\{\frac{p(p-1)}2^2, 0^{2p^2-2p}, \frac{p}2(-1-\sqrt5) ^{p-1},\frac{p}2(-1+\sqrt5) ^{p-1} \} 
& \mbox{ if } p\equiv 1  \\
\{\frac{p(p-1)}2^2, 0^{2p^2-2p}, \frac{p}2(-1-i\sqrt3) ^{p-1},\frac{p}2(-1+i\sqrt3) ^{p-1} \} 
& \mbox{ if } p\equiv 3  
\end{cases}\pmod{4}\\ 
\Lambda_5&=&\{p^1, 0^{2p-2}, -p^1, -\sqrt{p}^{p(p-1)}, \sqrt{p}^{p(p-1)}\}\\
\Lambda_6&=&\{p(p-1)^1, 0^{2p-2}, -p(p-1)^1, -\sqrt{p}^{p(p-1)}, \sqrt{p}^{p(p-1)}\}.
\end{eqnarray*}

As is customary, we use superscripts to indicate the multiplicity of each given eigenvalue. 

\medskip

\noindent{\bf Remark 5.}
It is worth mentioning that $\aut(\mathcal N_6)$ has constant rank 8. Due to this, 
all mergings of $\aut(\mathcal N_6)$ might be described with the aid of {\sf COCO}, 
at least for all values $p< 100$. Though this was accomplished by us for only a
few small values of $p$ (in parallel with computations in {\sf GAP}), such an approach  
might be helpful in conducting a more careful future analysis of similarly
obtained objects of constant rank.

\subsection{Rank 5 mergings}

It turns out that the Schurian rank 8 association scheme, which appears from the
2-orbits of $\aut(\mathcal N_6)$ has two non-Schurian rank 5 mergings. 
Relevant information for one of these schemes, denoted $\mathcal N_{5.1}$, is  
presented in the slides of M.K. Here, we shall restrict our attention only to 
consideration of the spectrum of $\mathcal N_{5.1}$. We believe it is of an 
independent interest due to its more sophisticated structure.

\medskip

\noindent{\bf Announcement 3.}
Denote again by $\Lambda_i$ the spectrum of the basic graph $\mathcal N_{5.1,i}$ in 
$\mathcal N_{5.1}$. Then
\begin{eqnarray*}
\Lambda_1 &=& \{ 1 ^{2p^2} \}\\
\Lambda_2 &=& \{  p-1^{2p}, -1^{2p^2-2p}\}\\
\Lambda_3 &=& \{ -p^{2p-2}, 0^{2p^2-2p}, p(p-1) ^2\}\\
\Lambda_4 &=& \left\{\pm\frac{p(p-1)}2^1, 0^{2p-2}, 
-\frac12\sqrt{p(p+1)}^{p(p-1)},\frac12\sqrt{p(p+1)}^{p(p-1)} \right\}\\
\Lambda_5 &=& \left\{\pm\frac{p(p+1)}2^1, 0^{2p-2},-\frac12\sqrt{p(p+1)}^{p(p-1)},
\frac12\sqrt{p(p+1)}^{p(p-1)} \right\}
\end{eqnarray*}

As before,  superscripts are used to indicate multiplicities of eigenvalues.

\subsection{Extra discussion}

Use of the term ``announcement''  in this section was a conscious decision 
on our part, intended to  stress the following:
\begin{itemize}
\item In the beginning, all results were obtained through plausible reasonings  
based on a careful analysis of numerous computer algebra experiments;
\item We are aware of all necessary tools in order to transform the results of plausible
reasonings to rigourously justified theoretical propositions;
\item Some such justifications have  already been reached, while others  
still require further effort and are postponed to the future;
\item A full justification of the formulated results does not fit into the agenda of
this paper, as much more will be needed in the way of  preliminaries and discussions 
of used techniques;
\item We expect to revisit this topic in the future and to devote to it a separate new paper;
\item This section intends also to underline our priority, at least in the formulation of
all the communicated results.
\end{itemize}

One of the essential features of the planned continuation of the presented research is 
working on the edge between AGT and extremal graph theory (briefly EGT). 
The established methodology for  determining  spectra of basic graphs is fairly 
traditional in the theory of association schemes,  while tools exploited in EGT 
are of a quite different nature. This is why we hope in our future work to present 
some fresh vision of a few classes of graphs exploited in EGT. There are also 
expectations that more careful analysis of some graphs in our schemes, 
in particular those of constant rank, may imply innovative results in the sphere of EGT.

\section{Concluding discussion}

\subsection{Origins of the project}

This research started in the framework of postdoctoral studies of the 
author \v S.Gy. at the Ben-Gurion University of the Negev, beginning in November 2011.
During the initial stages, \v S.Gy.\ was introduced to the main concepts of AGT
and was garnering new experience in  the use of a number of computer packages, 
as described  in Section 3. 

In particular, it was suggested to attempt to understand  without the use of a computer, 
the structure of two non-Schurian association schemes
on 18 points, which had been known for a long time.

Quite soon such a computer-free interpretation (strictly in terms of \cite{kr}) was 
elaborated. Moreover, it became clear that there was space for clever
generalisations. Further natural generalisations were elaborated, leading 
eventually to an understanding of the first model of the coherent 
configuration $\mathcal M$ and its four mergings.

At the next stage, ideas developed by F.\ Lazebnik, V.\ Ustimenko and A.J. Woldar, 
with which the author M.K.\ was already acquainted, were successfully exploited 
to understand the potential of the second model.
Moreover, both authors were sharing mutual pleasure from the ongoing feeling   
that the techniques of algebraic mergings outlined in \cite{km} works
perfectly well in our case also. 

Though this research was arranged without specific use of the techniques  
discussed in \cite{wo}, we believe it is quite likely that there is a very
promising potential impact from the amalgamation of the two approaches.

\subsection{More about some other computations}

While we think we have already devoted sufficient attention to diverse computer 
aided activities in previous sections, we nonetheless  feel compelled to at least 
briefly discuss a couple of extra approaches not touched upon earlier.
 
All  computations considered in these subsections were executed with the aid of
{\sf COCO~II}.

\subsubsection{Algebraic stabilizers}

The group $\aaut(\mathcal W)$ is a permutation group acting on the set of colors
of $\mathcal W$. Each merging corresponds to a partition $P$ of the set of colors
of $\mathcal W$. In effect, it is a set of sets of colors. 
If~$P$ is an algebraic merging then it consists of the orbits of some 
subgroup of $\aaut(\mathcal W)$. 
We can ask what is the \emph{algebraic stabilizer} of a partition $P$, 
i.e. the largest subgroup of $\aaut(\mathcal W)$ which leaves $P$ invariant.

The algebraic stabilizers of association schemes $\mathcal M_1, \mathcal M_2,
\mathcal M_3$ and $\mathcal M_4$ are  isomorphic to $\mathbb Z_2$, $\mathbb 
Z_2\times\mathbb Z_2$, $\mathbb Z_{p-1}\times \mathbb Z_2$ and $\mathbb 
Z_{p-1}\times D_4$, respectively.  (Here $D_4$ is the dihedral group of order $8$.) 
 
\subsubsection{Coherent subalgebras}

We computed numbers of coherent subalgebras of $\mathcal M$ which 
arise as algebraic mergings  for $p=3,5$ and $7$. 
We also attempted this for $p=11$, but while isomorphism testing for certain 
coherent subalgebras required only 1-2 minutes, others had to be  aborted, 
unsuccessfully, after a two-week period.
We  therefore limit our table below to numbers of pairwise non-isomorphic 
coherent configurations obtained as algebraic mergings of $\mathcal M$ for $p=3,5$ and $7$.

\medskip

\begin{center}
\begin{tabular}{|l|cccccc|}
\hline 
$p$ & CC & NCC & AS & Schur & NonSch & Intr \\
\hline
$p=3$ & 22 & 12 & 10 & 8 & 2 & 0 \\
$p=5$ & 60 & 36 & 24 & 18 & 6 & 3 \\
$p=7$ & 120 & 80 & 40 & 28 & 12 & 4 \\
\hline
\end{tabular}
\end{center}

\noindent We provide a legend for this table as follows: 
CC = number of coherent configurations, NCC = number of non-homogeneous
coherent configurations, AS = number of association schemes, Schur = 
number of Schurian association schemes, NonSch = number of non-Schurian
association schemes, Intr = number of non-Schurian association schemes
with intransitive group of automorphisms. 

We also arranged a number of computational experiments in order to enumerate
all homogeneous mergings of $\mathcal M$, again appearing as algebraic mergings,  
for a few small values of $p$. The results will be reported elsewhere. 

\subsection{Proofs versus plausible reasonings} 

Plausible reasonings (exactly in the sense of G.\ P\'olya \cite{po,pl}) 
play a very significant role in AGT in general, and in this project concretely. 
The modern use of computers, of course, extends the possibilities of this method 
of mathematical thinking. 

For decades, a striking model of research in AGT is to construct a new object with the
aid of a computer, to understand the object's properties (ideally to
create an aesthetically pleasing computer-free interpretation), and finally to extend the
initial object to an infinite series (or even class) of new structures. 

Fortunately, the authors  are able to report  exactly this level
of success in this paper. Of~course, {\em P\'olya's method}  (if  
fulfilled to its complete extent) ultimately requires rigourous proofs of all
claims, thus substituting plausible insights by justified
mathematical propositions. 

The authors are definitely on track to do this. Up to Section 11, we reported  
rigourous results, confirmed at the level of accepted standards of mathematical rigour. 
Completion of the wider project and its exposition  is simply a matter   
of time and space limitations.  We hope that the reader is able to detect, 
and even appreciate, our intentions.   

\subsection{About the style of this paper}

Our first and foremost task in this text is clear and traditional: to communicate
reached new results. While the authors of the majority of research papers are
typically satisfied by fulfilment of this goal, we were, from the very beginning,
thinking about some additional objectives   of our paper. 

Taking into account that both {\sf arXiv} and final versions are intended 
for online publication, space limitations are not so significant in such a case, 
although the reader's patience to read a longer text should definitely be taken 
under consideration. Striving to maintain a delicate balance between these two extrema, 
we briefly touch upon additional functions of communication of our text, 
formulated below as theses.

\subsubsection{On the edge with philosophy of mathematics}

\medskip

\noindent{\bf Education of a scientific researcher.}
As previously mentioned, the younger coauthor started this project with a quite modest goal:
to better explain known mathematical objects on 18 points. Fortunately for us, this
explanation was quickly extended to an understanding of new larger objects, 
ultimately leading to formal theoretical generalisations. 
Both authors shared the enjoyment of this efficient maturation of ideas, 
and firmly believe that their  process carries some unusual methodological
features. Sharing of them might be helpful for other mathematical parties.

\medskip

\noindent{\bf Philosophy of discovery.}
The computer was our closest ally at each stage of the  conducted research project. 
All software tools at our disposal were available free of charge. Moreover, 
when it became absolutely   necessary, we were able to communicate freely with 
colleagues who are high experts in the use of this software and benefit greatly 
from their advice. Although we could potentially conduct any number of experiments, 
usually 5-10 sufficed in order to guess a valid generalisation. At each stage it 
was a pleasure to experience  a gradual deepening of our understanding of new occurrences. 

\medskip

\noindent{\bf Role of proofs.} As a rule, we gained the main ingredients of our 
knowledge at a heuristic level. Typically, the elder coauthor was posing some 
quite natural questions in response to which the younger author was making 
the required calculations. A subsequent comparison of results would lead to 
the formulation of viable conjectures. With each successive stage of this 
methodological process, the younger colleague would  come to acquire a deeper 
appreciation and comprehension of this kind of mathematical behaviour. 
Proofs were  proceeded only after the corresponding picture became completely  
transparent to us, thus allowing the clearest possible formulation of  claims, 
already evident to us. 

\medskip

\noindent{\bf Advantages of a computer.} 
Once more, we were benefitting from the use of necessary free software. 
Typically, this was {\sf GAP} with its share packages, {\sf COCO}, {\sf COCO-IIR}, 
and programs from the home page of Hanaki and Miyamoto. 

We were also conducting experiments on how the resulting data should be presented 
and formatted, how tables could be organized for maximum benefit, making guesses of 
the structure of a group from its order, confirming  assumptions about groups with the 
aid of generators, etc.

\medskip

\noindent{\bf Goals of mathematics.}
This issue is on the edge with philosophy of science. One starts a very concrete project
within the well established area of AGT, applies some  traditional tools, 
and quickly becomes aware of the need to extend these tools in an innovative manner.
The obtained results, which are very often quite far from would have been 
predicted from the outset, leads one to investigate  
links between the exploited area and other parts of mathematics. 

What then should be regarded as the most significant result?
Is it discovery, proof, beauty of the results, new research horizons? 
It is difficult to say. Probably a combination of all  of these. 
This is how we naturally come to the forefront of philosophy of mathematics. 

\subsubsection{Extra vision from inside of philosophy of science}

Not wishing to open this box of Pandora, we just mention here a handful of 
interesting references and excerpts, to which we add our own brief personal reflections.

Its purpose is to guide the reader to a better understanding of our style 
of exposition and our three main methodological activities: computer 
experiments, subsequent generalisations, and proofs that we often found tedious.

\begin{itemize}
\item The paper \cite{th} played a very significant role in the long-standing discussion
between mathematicians and philosophers, which is here very briefly touched. It contains
thoughts of the late Fields medalist, with special attention to a ``continuing desire for human
understanding of a proof'' and advocating to make them as clear and simple as possible.
\item An interesting survey about the challenge of computer mathematics \cite{bwi} culminated
with a section entitled ``Romantic versus cool mathematics''.  Here \emph{romantic} refers to  
proofs that are navigable by mind, while  \emph{cool}  oppositely refers to proofs that are 
verified only by computer. 
\item The paper \cite{cc} is highly provocative from the outset, due  to its Section 1 title: 
``To prove or not to prove -- that is the question!''. 
The authors suggest to classify the evolution of a proof from ancient time to our millenium
in eight stages, especially attributing to the eighth and final stage a much more 
significant role of empirical and experimental features.  
\item Paper \cite{hn} was written by a famous expert in didactics of mathematics, and is 
devoted to the role of formal proof in high school education. Its presented 
list of functions of proof, which includes verification, explanation, discovery, communication,
etc., is fairly in concert with our own experience and philosophy, as they are discussed here.
\item Our last item refers to an article  \cite{wl} at the blog of a well known expert 
in computer algebra, and the creator and commercial promoter of {\sf Mathematica}. 
It contains an interesting discussion of knowledge obtained 
via computer, and speculates about the future of the role of proof in pure mathematics. 
\end{itemize}

\subsection{Bos\'ak graph}

There is one more fairly well known graph on 18 vertices, which may be easily 
obtained inside of our master coherent configuration on 18 points. 
This graph belongs to the family of directed strongly regular graphs, briefly DSRGs, 
a natural  generalisation of strongly regular (undirected) graphs to the case 
of mixed graphs. This concept was introduced by A.\ Duval in his seminal paper \cite{du}. 
We will use the established notation $(n,k,t,\lambda,\mu)$ for its parameter set.

The discussed graphs are regular graphs of valency $k$, and satisfy $AJ=JA=kJ$ and
$A^2=t\cdot I+\lambda\cdot A+\mu(J-I-A)$. Clearly, always $0\leq t\leq k$. 
If $t=k$ then we are getting  the  usual strongly regular graphs, while   
the case $t=0$ corresponds to doubly regular tournaments. 
For $0<t<k$ the wording \emph{genuine DSRG} was suggested. 

The main ingredients of the theory of DSRGs were developed in \cite{du} by Duval, 
who also discovered a few infinite classes of such graphs and posed problems 
of existence and full enumeration of all DSRGs, up to isomorphism, 
with a given parameter set. In particular, Duval mentioned the existence of a DSRG 
with the parameters $(18,4,3,0,1)$, and as well constructed an infinite family 
of such graphs with parameters $(k^2+k,k,1,0,1)$, where $k\geq 2$. 
It seems that Duval was unaware of the fact that his concrete graph 
on 18 vertices had already been discovered 10 years earlier by J. Bos\'ak 
(see \cite{bk,bo,bs}),  who was looking for so-called \emph{mixed Moore graphs}. 
In modern terms, these are DSRGs with $\lambda=0$ and $\mu=1$, which 
appear as a natural generalisation of classical (undirected) Moore graphs. 
For this class of mixed graphs, Bos\'ak developed a theory quite similar to 
the more general one established by Duval for all DSRGs later on.

It seems that for a long time Bos\'ak's results remained  undetected  by the
majority of experts in AGT. Fortunately, the authors of \cite{ng} rediscovered 
Bos\'ak's publications and breathed new life into them. 
In particular, they proved uniqueness of DSRGs with parameters $(18,4,3,0,1)$,  
and replied to some open questions posed by Bos\'ak. The results of Bos\'ak 
were definitely not known to the authors of \cite{fp} and \cite{kz}. 
They, in fact, independently duplicated  some of Bos\'ak's constructions 
in terms of 2-designs,  paying special attention to the case $(18,4,3,0,1)$. 
Nowadays we refer to the graph with these parameters as the \emph{Bos\'ak graph}~$B_{18}$. 
This graph is one of the main heroes in \cite{kz}; see Example 7.2 of that paper. 
The depiction of $B_{18}$ provided there vividly shows that it has $K_{3,3}$ 
as a quotient graph. Also it was shown in \cite{kz} that $G=\aut(B_{18})$ 
has order 108, and an explicit  set of generators for~$G$ was given. 
It was additionally  observed  that~$G$ is a central extension of~$\mathbb Z_3$ 
with the aid of the subgroup of index~2 in $\aut(K_{3,3})$, 
consisting of even permutations. (Incidentally, the original depiction of~$B_{18}$ 
given by Bos\'ak in \cite{bs} carries very much the same flavour as the one 
given in \cite{kz}, although we fell that the latter may be regarded as 
a bit more aesthetically pleasing.) 

The graph $B_{18}$ also attracted the attention of L.\ J\o rgensen \cite{jo}, 
who gave a description of~$B_{18}$ as a Cayley graph over a suitable (non-Abelian) 
group of order~18. We became aware of the references to Bos\'ak from the 
preliminary version of \cite{jo}, kindly sent to us by the author. 
This led the current authors to prepare the draft~\cite{gk} in which~$B_{18}$
is considered in the framework of our master configuration~$\mathcal M$ on~18 points. 
In particular, we showed that the coherent closure of~$B_{18}$ is a certain rank~7 
Schurian association scheme which arises as a merging of~$\mathcal M$.
The arc set of~$B_{18}$ is described as a union of  relations of~$\mathcal M$, 
specifically $A_1$, $C_2$, $E_0$ and $F_0$.  
Removal of all directed edges from~$B_{18}$ leads to the Pappus graph. 
Also a new geometric image of~$B_{18}$ was created  
which reflects the embedding of the Bos\'ak graph into the  torus. 
One more interesting finding was a graphical representation of~$B_{18}$ 
with the aid of the voltage graph of order~2 in the group $\mathbb Z_3\times\mathbb Z_3$.  

Last but not least, we wish to mention that the draft \cite{gk} provided 
the author \v S.Gy.\ a rare opportunity to enter into an additional research area of AGT.  
The recent joint publication~\cite{gy} provides evidence, as well as  hope, 
that his acquaintance with DSRGs is not a one-time affair. 

\subsection{Biaffine planes}

The authors learned the term ``biaffine plane'' from the paper \cite{wi} of P.\ Wild, 
who used this terminology in  a few of his earlier publications as well as in his 
Ph.D. thesis (University of~London, 1980).
Recently, however, we came to understand that the term was used quite earlier, 
namely in \cite{pi} by G.\ Pickert, one of the classic experts 
of modern finite geometries, who passed away in 2015 at the ripe old age of 97. 
The term is referred to M.\ Oehler (1975) and is discussed in the context of strongly regular
graphs, namely the famous Shrikhande's pseudo-$L_2$ association scheme.

A more extensive bibliographical search resulted in the discovery of a paper 
by~A.~Bennett~\cite{be} published in 1925. 
Quite in the style of that time the text \cite{be} does not
contain any references. Thus we urge experts in  the field of  history of finite geometries
to determine the earliest origins of this terminology.

\subsection{Spectra of our schemes}

In Appendix 3, we provide a description of the spectrum of each of the four 
association schemes $\mathcal M_1$--$\mathcal M_4$ determined in this paper. 
As was discussed earlier, a justification of this result (as with similar results 
for rank~6 and rank~5 schemes) does not fit well the framework of this text.  
Nevertheless, we hope that by affording this spectral information to others, 
it may  be of some definite help. We alert the reader that 
the last row in some tables still lacks a precise formulation.

\subsection{A few extra references}

We attract the reader's attention to a few  recent publications which provide some
interesting overlap with the topics touched here by us.

The paper \cite{ll} belongs to the area of coding theory. 
The introduced LDPC codes are based on an infinite family of bipartite graphs, 
known notationally as $D(2,q)$, which initially arose in the context 
of extremal graph theory (EGT). In \cite{ll} the spectrum of $D(2,q)$ is determined, 
and it is proved that each such graph is Ramanujan. However, the reported spectrum 
is in error as it is not symmetric about 0, a theoretical requirement for bipartite graphs. 
Clearly, for $q$ prime the graphs $D(2,q)$ are living inside our scheme $\mathcal M(q)$.
Thus it would be of special interest to compare the arguments in~\cite{ll} 
with the data provided in our text.
  
The paper \cite{br} deals with graphs related to cages. Here again the 
considered graphs are living inside our configuration~$\mathcal M_p$, and it 
is shown that they are close to optimal in the framework of~EGT.

Authors of the paper \cite{ci} describe the spectrum of the Wenger graphs $W_m(q)$. 
Some formulas for the multiplicities of the eigenvalues are provided. For certain 
cases, a comparison of these values to the ones that result from our own formulas 
might be of ample curiosity.

\subsection{Presentations} Partial reports about the results expressed in this text were
presented a few times:
\begin{itemize}
\item by M.K.\ at  the workshop ``84th workshop on general algebra'', Dresden, Germany, 
June 2012, see {\tt http://tu-dresden.de/die\_tu\_dresden/fakultaeten/
\newline 
fakultaet\_mathematik\_und\_naturwissenschaften/fachrichtung\_mathematik/
\newline
institute/algebra/aaa84/};
\item by \v S.Gy.\ at the conference ``Computers in Scientific Discovery 6'', 
Portoro\v z, Slovenia, 2012, see 
{\tt http://conferences2.imfm.si/internalPage.py?pageId=15\&confId=12};
\item by \v S.Gy.\ at the workshop ``50th Summer school on general algebra and ordered sets'',
Nov\'y Smokovec, Slovakia, Sept.\ 2012, see {\tt https://sites.google.com/site/ssalgebra2012/};
\item by M.K.\ at the seminar of Queen Mary, University of London, England, May 2013, 
see {\tt http://www.maths.qmul.ac.uk/seminars/};
\item by M.K.\ at the conference ``Modern trends in AGT'', Villanova (PA), USA, June 2014, 
see {\tt https://www1.villanova.edu/villanova/artsci/mathematics/ \newline
newsevents/mtagt/slides-of-all-talks.html}.
\end{itemize}

Each of these presentations was very helpful for the authors, giving them the chance to better
understand ways in which they could improve successive versions of the paper.

\subsection{Research in progress}

As was mentioned, we see a very promising intersection of our research methodology 
with activities in extremal graph theory. We are currently in a position to exploit 
this potential to its full extent. In particular, we aim to investigate how families 
of association schemes discovered by us (as well as other possible examples and classes of
such structures) may be ``repurposed'' for beneficial application to other diverse 
branches of graph theory.

\section*{Acknowledgements}

The stay of the first author at the Ben-Gurion University 
and cooperation between the authors was enabled and supported by 
a scholarship based on an agreement between the Israeli and Slovak governments. 
The first author gratefully acknowledges the contribution 
of the Scientific Grant Agency of the Slovak Republic under the grant 1/1005/12. 
The second author was partially supported by the EUROCORES Programme 
EUROGIGA (project GReGAS) of the European Science Foundation.

Both authors were supported by the Project: Mobility - enhancing research, 
science and education at the Matej Bel University, ITMS code: 26110230082, 
under the Operational Program Education cofinanced by the European Social Fund.

The author M.K. is pleased to thank his friends and colleagues Felix 
Lazebnik, Vasyl Ustimenko, and especially Andrew Woldar for striking 
long-standing discussions related to their methods and results.  In addition, 
both authors are grateful to Woldar for his careful linguistic polishing 
of earlier versions  of the text. 

We thank Misha Muzychuk for helpful conversations. Cooperation in computer
algebra provided by Danny Kalmanovich, Christian Pech, Sven Reichard and Matan
Ziv-Av is gratefully acknowledged. Finally, we thank Jozef \v Sir\'a\v n and 
M\'aria \v Zd\'\i malov\'a for helpful communication. 

\section*{Appendix 1}

In this appendix we would like to display the intersection numbers of the
mentioned association schemes. 
For the sake of brevity let us denote 
$\xi:=\delta_{i+j,k}+\delta_{i-j,k}+\delta_{-i+j,k}+\delta_{-i-j,k}$,
and $M_{ijk}:=\max\{\delta_{i+j,k},\delta_{i-j,k}\}+\max\{\delta_{-i+j,k},\delta_{-i-j,k}\}$.

The superscript in the upper left corner indicates which relation is fixed in the table.
The subscripts $i$ and $j$ correspond to rows and columns, respectively, of the 
intersection matrices. 

\medskip

\noindent{\bf Intersection numbers of $\mathcal M_1$:}

\medskip

\begin{minipage}{7.5cm}
\begin{center}
\begin{tabular}{|c|cccc|}
\hline
$c_{i,j}^0$ & $0$ & $S_j$ & $T_j$ & $U_j$ \\
\hline
$0$ & $1$ & $0$ & $0$ & $0$ \\
$S_i$ & $0$ & $\delta_{i,j'}$ & $0$ & $0$ \\
$T_i$ & $0$ & $0$ & $p\cdot\delta_{i,j'}$ & $0$ \\
$U_i$ & $0$ & $0$ & $0$ & $p\cdot\delta_{i,j'}$ \\
\hline
\end{tabular}
\end{center}
\end{minipage}
\begin{minipage}{7.5cm}
\begin{center}
\begin{tabular}{|c|cccc|}
\hline
$c_{i,j}^{sk}$ & $0$ & $S_j$ & $T_j$ & $U_j$ \\
\hline
$0$ & $0$ & $\delta_{j,k}$ & $0$ & $0$ \\
$S_i$ & $\delta_{i,k}$ & $\delta_{i+j,k}$ & $0$ & $0$ \\
$T_i$ & $0$ & $0$ & $p\cdot\delta_{i,j'}$ & $0$ \\
$U_i$ & $0$ & $0$ & $0$ & $p\cdot\delta_{i+j,k}$ \\
\hline
\end{tabular}
\end{center}
\end{minipage}

\medskip

\begin{minipage}{7.5cm}
\begin{center}
\begin{tabular}{|c|cccc|}
\hline
$c_{i,j}^{tk}$ & $0$ & $S_j$ & $T_j$ & $U_j$ \\
\hline
$0$ & $0$ & $0$ & $\delta_{j,k}$ & $0$ \\
$S_i$ & $0$ & $0$ & $\delta_{j,k}$ & $0$ \\
$T_i$ & $\delta_{i,k}$ & $\delta_{i,k}$ & $p\cdot\delta_{i+j,k}$ & $0$ \\
$U_i$ & $0$ & $0$ & $0$ & $1$ \\
\hline
\end{tabular}
\end{center}
\end{minipage}
\begin{minipage}{7.5cm}
\begin{center}
\begin{tabular}{|c|cccc|}
\hline
$c_{i,j}^{uk}$ & $0$ & $S_j$ & $T_j$ & $U_j$ \\
\hline
$0$ & $0$ & $0$ & $0$ & $\delta_{j,k}$ \\
$S_i$ & $0$ & $0$ & $0$ & $\delta_{i+j,k}$ \\
$T_i$ & $0$ & $0$ & $0$ & $1$ \\
$U_i$ & $\delta_{i,k}$ & $\delta_{i+j,k}$ & $1$ & $0$ \\
\hline
\end{tabular}
\end{center}
\end{minipage}

\medskip

\noindent{\bf Intersection numbers of $\mathcal M_2$:}

\medskip

\begin{minipage}{7.5cm}
\begin{center}
\begin{tabular}{|c|ccccc|}
\hline
$c_{i,j}^0$ & $0$ & $S^*_j$ & $T_j$ & $U_0$ & $U^*_j$ \\
\hline
$0$ & $1$ & $0$ & $0$ & $0$ & $0$ \\
$S^*_i$ & $0$ & $2\cdot\delta_{i,j}$ & $0$ & $0$ & $0$ \\
$T_i$ & $0$ & $0$ & $p\cdot\delta_{i,j'}$ & $0$ & $0$ \\
$U_0$ & $0$ & $0$ & $0$ & $p$ & $0$ \\
$U^*_i$ & $0$ & $0$ & $0$ & $0$ & $2p\cdot\delta_{i,j}$ \\
\hline
\end{tabular}
\end{center}
\end{minipage}
\begin{minipage}{7.5cm}
\begin{center}
\begin{tabular}{|c|ccccc|}
\hline
$c_{i,j}^{sk}$ & $0$ & $S^*_j$ & $T_j$ & $U_0$ & $U^*_j$ \\
\hline
$0$ & $0$ & $\delta_{j,k}$ & $0$ & $0$ & $0$ \\
$S^*_i$ & $\delta_{i,k}$ & $\xi$ & $0$ & $0$ & $0$ \\
$T_i$ & $0$ & $0$ & $p\cdot\delta_{i,j'}$ & $0$ & $0$ \\
$U_0$ & $0$ & $0$ & $0$ & $p\cdot\xi$ & $p\cdot\xi$ \\
$U^*_i$ & $0$ & $0$ & $0$ & $p\cdot\xi$ & $p\cdot\xi$ \\
\hline
\end{tabular}
\end{center}
\end{minipage}

\medskip

\begin{minipage}{7.5cm}
\begin{center}
\begin{tabular}{|c|ccccc|}
\hline
$c_{i,j}^{tk}$ & $0$ & $S^*_j$ & $T_j$ & $U_0$ & $U^*_j$ \\
\hline
$0$ & $0$ & $0$ & $\delta_{j,k}$ & $0$ & $0$ \\
$S^*_i$ & $0$ & $0$ & $2\cdot \delta_{j,k}$ & $0$ & $0$ \\
$T_i$ & $\delta_{i,k}$ & $2\cdot \delta_{i,k}$ & $p\cdot\delta_{i+j,k}$ & $0$ & $0$ \\
$U_0$ & $0$ & $0$ & $0$ & $1$ & $2$ \\
$U^*_i$ & $0$ & $0$ & $0$ & $2$ & $4$ \\
\hline
\end{tabular}
\end{center}
\end{minipage}
\begin{minipage}{7.5cm}
\begin{center}
\begin{tabular}{|c|ccccc|}
\hline
$c_{i,j}^{uk}$ & $0$ & $S^*_j$ & $T_j$ & $U_0$ & $U^*_j$ \\
\hline
$0$ & $0$ & $0$ & $0$ & $\delta_{0,k}$ & $\delta_{0,j}$ \\
$S^*_i$ & $0$ & $0$ & $0$ & $M_{i0k}$ & $M_{ijk}$ \\
$T_i$ & $0$ & $0$ & $0$ & $1$ & $2$ \\
$U_0$ & $\delta_{0,k}$ & $M_{0jk}$ & $1$ & $0$ & $0$ \\
$U^*_i$ & $\delta_{0,i}$ & $M_{ijk}$ & $2$ & $0$ & $0$ \\
\hline
\end{tabular}
\end{center}
\end{minipage}

\medskip

\noindent{\bf Intersection numbers of $\mathcal M_3$:}

\medskip

\begin{minipage}{7.5cm}
\begin{center}
\begin{tabular}{|c|ccccc|}
\hline
$c_{i,j}^0$ & $0$ & $S$ & $T_j$ & $U_0$ & $U$ \\
\hline
$0$ & $1$ & $0$ & $0$ & $0$ & $0$ \\
$S$ & $0$ & $p-1$ & $0$ & $0$ & $0$ \\
$T_i$ & $0$ & $0$ & $p\cdot\delta_{i,j'}$ & $0$ & $0$ \\
$U_0$ & $0$ & $0$ & $0$ & $p$ & $0$ \\
$U$ & $0$ & $0$ & $0$ & $0$ & $p(p-1)$ \\
\hline
\end{tabular}
\end{center}
\end{minipage}
\begin{minipage}{7.5cm}
\begin{center}
\begin{tabular}{|c|ccccc|}
\hline
$c_{i,j}^{s}$ & $0$ & $S$ & $T_j$ & $U_0$ & $U$ \\
\hline
$0$ & $0$ & $1$ & $0$ & $0$ & $0$ \\
$S$ & $1$ & $p-2$ & $0$ & $0$ & $0$ \\
$T_i$ & $0$ & $0$ & $p\cdot\delta_{i,j'}$ & $0$ & $0$ \\
$U_0$ & $0$ & $0$ & $0$ & $0$ & $p$ \\
$U$ & $0$ & $0$ & $0$ & $p$ & $p(p-2)$ \\
\hline
\end{tabular}
\end{center}
\end{minipage}

\medskip

\begin{minipage}{15cm}
\begin{center}
\begin{tabular}{|c|ccccc|}
\hline
$c_{i,j}^{tk}$ & $0$ & $S$ & $T_j$ & $U_0$ & $U$ \\
\hline
$0$ & $0$ & $0$ & $\delta_{j,k}$ & $0$ & $0$ \\
$S$ & $0$ & $0$ & $(p-1)\cdot \delta_{j,k}$ & $0$ & $0$ \\
$T_i$ & $\delta_{i,k}$ & $(p-1)\cdot \delta_{i,k}$ & $p\cdot\delta_{i+j,k}$ & $0$ & $0$ \\
$U_0$ & $0$ & $0$ & $0$ & $1$ & $p-1$ \\
$U$ & $0$ & $0$ & $0$ & $p-1$ & $(p-1)^2$ \\
\hline
\end{tabular}
\end{center}
\end{minipage}

\medskip

\begin{minipage}{7.5cm}
\begin{center}
\begin{tabular}{|c|ccccc|}
\hline
$c_{i,j}^{u0}$ & $0$ & $S$ & $T_j$ & $U_0$ & $U$ \\
\hline
$0$ & $0$ & $0$ & $0$ & $1$ & $0$ \\
$S$ & $0$ & $0$ & $0$ & $0$ & $p-1$ \\
$T_i$ & $0$ & $0$ & $0$ & $1$ & $p-1$ \\
$U_0$ & $1$ & $0$ & $1$ & $0$ & $0$ \\
$U$ & $0$ & $p-1$ & $p-1$ & $0$ & $0$ \\
\hline
\end{tabular}
\end{center}
\end{minipage}
\begin{minipage}{7.5cm}
\begin{center}
\begin{tabular}{|c|ccccc|}
\hline
$c_{i,j}^{u}$ & $0$ & $S$ & $T_j$ & $U_0$ & $U$ \\
\hline
$0$ & $0$ & $0$ & $0$ & $0$ & $1$ \\
$S$ & $0$ & $0$ & $0$ & $1$ & $p-2$ \\
$T_i$ & $0$ & $0$ & $0$ & $1$ & $p-1$ \\
$U_0$ & $0$ & $1$ & $1$ & $0$ & $0$ \\
$U$ & $1$ & $p-2$ & $p-1$ & $0$ & $0$ \\
\hline
\end{tabular}
\end{center}
\end{minipage}

\newpage

\noindent{\bf Intersection numbers of $\mathcal M_4$:}

\medskip

\begin{minipage}{7.5cm}
\begin{center}
\begin{tabular}{|c|ccccc|}
\hline
$c_{i,j}^0$ & $0$ & $S$ & $T^*_j$ & $U_0$ & $U$ \\
\hline
$0$ & $1$ & $0$ & $0$ & $0$ & $0$ \\
$S$ & $0$ & $p-1$ & $0$ & $0$ & $0$ \\
$T^*_i$ & $0$ & $0$ & $2p\cdot\delta_{i,j}$ & $0$ & $0$ \\
$U_0$ & $0$ & $0$ & $0$ & $p$ & $0$ \\
$U$ & $0$ & $0$ & $0$ & $0$ & $p(p-1)$ \\
\hline
\end{tabular}
\end{center}
\end{minipage}
\begin{minipage}{7.5cm}
\begin{center}
\begin{tabular}{|c|ccccc|}
\hline
$c_{i,j}^{s}$ & $0$ & $S$ & $T^*_j$ & $U_0$ & $U$ \\
\hline
$0$ & $0$ & $1$ & $0$ & $0$ & $0$ \\
$S$ & $1$ & $p-2$ & $0$ & $0$ & $0$ \\
$T^*_i$ & $0$ & $0$ & $2p\cdot\delta_{i,j}$ & $0$ & $0$ \\
$U_0$ & $0$ & $0$ & $0$ & $0$ & $p$ \\
$U$ & $0$ & $0$ & $0$ & $p$ & $p(p-2)$ \\
\hline
\end{tabular}
\end{center}
\end{minipage}

\medskip

\begin{minipage}{15cm}
\begin{center}
\begin{tabular}{|c|ccccc|}
\hline
$c_{i,j}^{tk}$ & $0$ & $S$ & $T^*_j$ & $U_0$ & $U$ \\
\hline
$0$ & $0$ & $0$ & $\delta_{j,k}$ & $0$ & $0$ \\
$S$ & $0$ & $0$ & $(p-1)\cdot \delta_{j,k}$ & $0$ & $0$ \\
$T^*_i$ & $\delta_{i,k}$ & $(p-1)\cdot \delta_{i,k}$ & $p\cdot\xi$ & $0$ & $0$ \\
$U_0$ & $0$ & $0$ & $0$ & $1$ & $p-1$ \\
$U$ & $0$ & $0$ & $0$ & $p-1$ & $(p-1)^2$ \\
\hline
\end{tabular}
\end{center}
\end{minipage}

\medskip

\begin{minipage}{7.5cm}
\begin{center}
\begin{tabular}{|c|ccccc|}
\hline
$c_{i,j}^{u0}$ & $0$ & $S$ & $T^*_j$ & $U_0$ & $U$ \\
\hline
$0$ & $0$ & $0$ & $0$ & $1$ & $0$ \\
$S$ & $0$ & $0$ & $0$ & $0$ & $p-1$ \\
$T_i$ & $0$ & $0$ & $0$ & $2$ & $2p-2$ \\
$U_0$ & $1$ & $0$ & $2$ & $0$ & $0$ \\
$U$ & $0$ & $p-1$ & $2p-2$ & $0$ & $0$ \\
\hline
\end{tabular}
\end{center}
\end{minipage}
\begin{minipage}{7.5cm}
\begin{center}
\begin{tabular}{|c|ccccc|}
\hline
$c_{i,j}^{u}$ & $0$ & $S$ & $T^*_j$ & $U_0$ & $U$ \\
\hline
$0$ & $0$ & $0$ & $0$ & $0$ & $1$ \\
$S$ & $0$ & $0$ & $0$ & $1$ & $p-2$ \\
$T^*_i$ & $0$ & $0$ & $0$ & $2$ & $2p-2$ \\
$U_0$ & $0$ & $1$ & $2$ & $0$ & $0$ \\
$U$ & $1$ & $p-2$ & $2p-2$ & $0$ & $0$ \\
\hline
\end{tabular}
\end{center}
\end{minipage}

\section*{Appendix 2}

The table below reflects the manner in which we were able to recognize the 
existence of the non-Schurian association scheme $\mathcal N_6$ based on  
analysis of computer data. The column headings indicate the value of $p$, 
the order of the scheme, the order of the automorphism group of the Schurian rank 8 
scheme, the order of $\aut(\mathcal N_6)$, and whether or not $\mathcal N_6$ is
commutative, or symmetric.
 
\begin{center}
\begin{tabular}{|r|r|r|r|c|c|}
\hline
$p$ & $n$ & rank 8 & rank 6 & commutative & symmetric \\
\hline
3 & 18 & 54 & 54 & yes & no \\
5 & 50 & 1,000 & 1,000 & yes & yes \\
7 & 98 & 6,174 & 6,174 & yes & no \\
11 & 242 & 66,550 & 66,550 & yes & no \\
13 & 338 & 158,184 & 158,184 & yes & yes \\
17 & 578 & 628,864 & 628,864 & yes & yes \\
19 & 722 & 1,111,158 & 1,111,158 & yes & no \\
\hline
\end{tabular}
\end{center}

\section*{Appendix 3: Spectrum of schemes $\mathcal M_1$--$\mathcal M_4$}

In all association schemes the spectrum of the identity matrix of order $2p^2$ is omitted. 

\subsection*{Association scheme $\mathcal M_1$}

There are $p-1$ basic matrices with spectrum 

\begin{center}
\begin{tabular}{c|l}
eigenvalue & multiplicity \\
\hline
1 & $2p$ \\
roots of $1+x+x^2+\ldots+x^{p-1}$ & $2p$ \\ 
\end{tabular}
\end{center}

$p-1$ matrices with spectrum

\begin{center}
\begin{tabular}{c|l}
eigenvalue & multiplicity \\
\hline
$p$ & 2 \\
0 & $2p(p-1)$ \\
roots of $\displaystyle\sum_{i=0}^{p-1}p^i\cdot x^{p-1-i}$ & 2 \\ 
\end{tabular}
\end{center}

one matrix with spectrum

\begin{center}
\begin{tabular}{c|l}
eigenvalue & multiplicity \\
\hline
$\!\!\!\! -p$ & 1 \\
$p$ & 1 \\
0 & $2p-2$ 
\\
$\!\!\!\! -\sqrt{p}$ & $p(p-1)$ \\
$\sqrt{p}$ & $p(p-1)$ \\
\end{tabular}
\end{center}

$p-1$ matrices with spectrum 

\begin{center}
\begin{tabular}{c|l}
eigenvalue & multiplicity \\
\hline
$\!\!\!\! -p$ & 1 \\
$p$ & 1 \\
0 & $2p-2$ \\
roots of some polynomial of degree $2p-2$ & $p$ \\ 
\end{tabular}
\end{center}

\subsection*{Association scheme $\mathcal M_2$}

There are $(p-1)/2$ matrices with spectrum

\begin{center}
\begin{tabular}{c|l}
eigenvalue & multiplicity \\
\hline
$p-1$ & $2p$ \\
roots of $x^{(p-1)/2}+x^{(p-3)/2}-\frac{p-3}2 x^{(p-5)/2}-\frac{p-5}2 x^{(p-7)/2}+\ldots$ & $4p$ \\ 
\end{tabular}
\end{center}

$p-1$ matrices with spectrum

\begin{center}
\begin{tabular}{c|l}
eigenvalue & multiplicity \\
\hline
$p$ & 2 \\
0 & $2p(p-1)$ \\
roots of $\displaystyle\sum_{i=0}^{p-1}p^i\cdot x^{p-1-i}$ & 2 \\ 
\end{tabular}
\end{center}

one matrix with spectrum

\begin{center}
\begin{tabular}{c|l}
eigenvalue & multiplicity \\
\hline
$\!\!\!\! -p$ & 1 \\
$p$ & 1 \\
0 & $2p-2$ 
\\
$\!\!\!\! -\sqrt{p}$ & $p(p-1)$ \\
$\sqrt{p}$ & $p(p-1)$ \\
\end{tabular}
\end{center}

$p-1$ matrices with spectrum 

\begin{center}
\begin{tabular}{c|l}
eigenvalue & multiplicity \\
\hline
$\!\!\!\! -2p$ & 1 \\
$2p$ & 1 \\
0 & $2p-2$ \\
roots of some polynomial of degree $p-1$ & $2p$ \\
\end{tabular}
\end{center}

\subsection*{Association scheme $\mathcal M_3$}

There is one matrix with spectrum

\begin{center}
\begin{tabular}{c|l}
eigenvalue & multiplicity \\
\hline
$p-1$ & $2p$ \\
$-1$ & $2p(p-1)$ \\ 
\end{tabular}
\end{center}

$p-1$ matrices with spectrum

\begin{center}
\begin{tabular}{c|l}
eigenvalue & multiplicity \\
\hline
$p$ & 2 \\
0 & $2p(p-1)$ \\
roots of $\displaystyle\sum_{i=0}^{p-1}p^i\cdot x^{p-1-i}$ & 2 \\ 
\end{tabular}
\end{center}

one matrix with spectrum

\begin{center}
\begin{tabular}{c|l}
eigenvalue & multiplicity \\
\hline
$\!\!\!\! -p$ & 1 \\
$p$ & 1 \\
0 & $2p-2$ \\
$\!\!\!\! -\sqrt{p}$ & $p(p-1)$ \\
$\sqrt{p}$ & $p(p-1)$ \\
\end{tabular}
\end{center}

one matrix with spectrum 

\begin{center}
\begin{tabular}{c|l}
eigenvalue & multiplicity \\
\hline
$\!\!\!\! -p(p-1)$ & 1 \\
$\!\!\!\! -\sqrt{p}$ & $p(p-1)$ \\
0 & $2p-2$ \\
$\sqrt{p}$ & $p(p-1)$ \\
$p(p-1)$ & 1 \\
\end{tabular}
\end{center}

\subsection*{Association scheme $\mathcal M_4$}

There is one matrix with spectrum

\begin{center}
\begin{tabular}{c|l}
eigenvalue & multiplicity \\
\hline
$p-1$ & $2p$ \\
$-1$ & $2p(p-1)$ \\ 
\end{tabular}
\end{center}

$(p-1)/2$ matrices with spectrum

\begin{center}
\begin{tabular}{c|l}
eigenvalue & multiplicity \\
\hline
$2p$ & 2 \\
0 & $2p(p-1)$ \\
roots of some polynomial of degree $(p-1)/2$ & 4 \\ 
\end{tabular}
\end{center}

one matrix with spectrum

\begin{center}
\begin{tabular}{c|l}
eigenvalue & multiplicity \\
\hline
$\!\!\!\! -p$ & 1 \\
$p$ & 1 \\
0 & $2p-2$ \\
$\!\!\!\! -\sqrt{p}$ & $p(p-1)$ \\
$\sqrt{p}$ & $p(p-1)$ \\
\end{tabular}
\end{center}

one matrix with spectrum 

\begin{center}
\begin{tabular}{c|l}
eigenvalue & multiplicity \\
\hline
$\!\!\!\! -p(p-1)$ & 1 \\
$\!\!\!\! -\sqrt{p}$ & $p(p-1)$ \\
0 & $2p-2$ \\
$\sqrt{p}$ & $p(p-1)$ \\
$p(p-1)$ & 1 \\
\end{tabular}
\end{center}

\end{document}